\documentclass[a4paper,
               10pt,
               pdftex,
               normalheadings,
               headsepline,
               footsepline,
               headinclude,
               footinclude,
               DIV14,
               abstracton]
{scrartcl}

\usepackage[left=2.50cm, right=2.50cm, top=2.00cm, bottom=3.50cm]{geometry}

\usepackage{authblk}
\usepackage[ngerman, english]{babel}

\usepackage{amsmath,nicefrac}
\usepackage{amssymb}
\usepackage{bm}
\usepackage{algorithmic}
\usepackage{algorithm} 
\usepackage{array}
\usepackage{mathtools}
\usepackage{enumerate}
\usepackage{xcolor}
\usepackage[english]{babel}
\usepackage{placeins}
\usepackage{caption}
\usepackage{subcaption}
\usepackage{multirow}
\usepackage{thmtools}
\usepackage{thm-restate}
\usepackage{array}
\usepackage{makecell}

\usepackage{color}
\usepackage{tikz}
\usepackage{cite}
\usepackage[normalem]{ulem} 

\newcommand{\WW}{\Xi}

\newcommand{\Hspace}{L^2(\Omega)}
\newcommand{\Vspace}{H^1_0(\Omega)}
\newcommand{\Vdual}{H^{-1}(\Omega)}

\definecolor{DarkOrange}{rgb}{0.83 0.33 0}

\DeclareMathOperator*{\argmin}{arg\,min}
\DeclareMathOperator*{\conv}{conv}
\newcommand{\closure}[1]{{\overline{#1}}}
\newcommand{\R}{\mathbb{R}}
\newcommand{\N}{\mathbb{N}}
\renewcommand{\H}{\mathcal{H}}
\newcommand{\dual}{{\ensuremath{^*}}}

\usepackage{amsthm}
\newtheorem{theorem}{Theorem}[section]
\newtheorem{definition}[theorem]{Definition}
\newtheorem{proposition}[theorem]{Proposition}
\newtheorem{lemma}[theorem]{Lemma}
\newtheorem{remark}[theorem]{Remark}
\newtheorem{corollary}[theorem]{Corollary}
\newtheorem{example}[theorem]{Example}

\begin{document}

\title{A Descent Method for Nonsmooth Multiobjective Optimization in Hilbert Spaces}
\author[1]{Konstantin Sonntag}
\author[1]{Bennet Gebken}
\author[2]{Georg Müller}
\author[3]{Sebastian Peitz}
\author[4]{Stefan Volkwein}
\affil[1]{\normalsize Department of Mathematics, Paderborn University, Germany}
\affil[2]{\normalsize Faculty of Mathematics, Heidelberg University, Germany}
\affil[3]{\normalsize Department of Computer Science, Paderborn University, Germany}
\affil[4]{\normalsize Department of Mathematics and Statistics, Konstanz University, Germany}

\date{}

\maketitle

\begin{abstract}
\noindent The efficient optimization method for locally Lipschitz continuous multiobjective optimization problems from \cite{G2021} is extended from finite-dimensional problems to general Hilbert spaces. The method iteratively computes Pareto critical points, where in each iteration, an approximation of the subdifferential is computed in an efficient manner and then used to compute a common descent direction for all objective functions. To prove convergence, we present some new optimality results for nonsmooth multiobjective optimization problems in Hilbert spaces. Using these, we can show that every accumulation point of the sequence generated by our algorithm is Pareto critical under common assumptions. Computational efficiency for finding Pareto critical points is numerically demonstrated for multiobjective optimal control of an obstacle problem.
\end{abstract}
\section{Introduction}

There are many applications, where multiple objectives have to be optimized at the same time. For example, when manufacturing a product, one wants to maximize the quality and simultaneously minimize the production cost. This leads to a \emph{multiobjective optimization problem} (MOP), where the goal is to find all optimal compromises between the objectives. Naturally, there are applications, where the objectives are defined on Hilbert spaces and feature nonsmoothness. For example, in \cite{MM1993}, an obstacle problem with an elastic string is considered, where one objective is maximization of the contact area between the string and a given obstacle and another objective is minimization of the total force applied to the string.

There is a vast amount of methods available for solving various types of \emph{finite}-dimensional optimization problems, but while most of them are designed to deal with \emph{either} nonsmoothness (e.g., \cite{BKM2014}) \emph{or} multiple objectives (e.g., \cite{Miettinen1998, E2021}), algorithms for nonsmooth MOPs are scarce. Two possible methods designed for nonsmooth MOPs are the \emph{proximal bundle method} \cite{MM2003, MKW2014} and the \emph{gradient sampling method} \cite{G2021}. 

Combining nonsmoothness, multiple objectives \emph{and} an infinite-dimensional Hilbert space setting becomes additionally challenging.
When presented with such a nonsmooth MOP in infinite-dimensions, there are several options to proceed, among them:
\begin{enumerate}
    \item \label{opt:disc} Discretize the infinite-dimensional nonsmooth MOP and then use a solver for finite dimensional problems, e.g., one of those presented in \cite{MM2003, MKW2014,G2021}.
    \item \label{opt:scalar} Scalarize the problem and then use a solver for infinite-dimensional non-smooth scalar optimization, e.g., \cite{Miettinen1998, Bernreuther2022}.
    \item \label{opt:awesome} Design a method that is capable of treating infinite dimensions, nonsmooth objective functions and multiple objectives at the same time.
\end{enumerate}
Option \ref{opt:disc} does not incorporate the underlying infinite-dimensional problem's topology and can therefore suffer from mesh-dependent behavior such as inconsistent termination criteria between different meshes; cf., e.g., the discussion in \cite[Sections~3.2.2-3.2.4]{HPUU09}. Option \ref{opt:scalar}, as in the smooth case, struggles in the presence of nonconvexity or when the number of objectives exceeds two. Option \ref{opt:awesome} suffers from neither of these drawbacks but is technically challenging to realize, and while infinite-dimensional nonsmooth MOPs with additional structure, such as convexity or composite structure, have previously been addressed (e.g., \cite{BIS2005, G2019, BG2018}), to the best of the authors' knowledge, there are no nonscalarizing methods for solving general, unstructured nonsmooth infinite-dimensional MOPs. 

The goal of this article is to generalize the common descent method based on subderivative sampling presented in \cite{G2021} from  finite-dimensional to infinite-dimensional (Hilbert space) settings. The main idea in \cite{G2021} is to replace the Clarke subdifferential \cite{C1983} in the design of the descent direction in the dynamic gradient approach of \cite{AGG2015} with the Goldstein $\varepsilon$-subdifferential \cite{G1977}, and to approximate the latter via an adaptive gradient sampling scheme. This way, a descent direction for nonsmooth MOPs can be computed. Combining this descent direction with an Armijo-backtracking-type step size control yields a descent method, for which convergence to points satisfying a necessary optimality condition has been shown. This algorithmic approach can be extended to a general Hilbert space setting in a relatively straight-forward manner, but the convergence analysis of the algorithm requires modifications to account for the loss of compactness. Additionally, the notions of optimality employed in \cite{G2021} will be adapted. While the Clarke subdifferential and the Goldstein $\varepsilon$-subdifferential have already been defined on Hilbert spaces \cite{C1983,M2018,M2019}, their multiobjective counterparts require additional attention. We generalize these objects and prove that they have a generalized demi-closedness property, and employ them in the derivation of necessary conditions for Pareto optimality.

This article is organized as follows. In Section~\ref{sec:Theoretical Background}, we introduce the basics of multiobjective optimization and nonsmooth analysis in Hilbert spaces. In Section~\ref{sec:eps-subgradient for MOPs}, we extend the Goldstein $\varepsilon$-subdifferential to the multiobjective, infinite-dimensional setting and investigate its properties. Theorem~\ref{thm:subdiff_eps_subdiff_multi} describes a demi-closedness property of the multiobjective $\varepsilon$-subdifferential, which is important for the convergence proof of the introduced method. The main results of this article are presented in Section~\ref{sec:Descent_method}. First, we describe how descent directions satisfying a sufficient descent property for all objective functions can be obtained theoretically using the extended subdifferential from the previous section. In Subsection~\ref{subsec:Efficient Computation of Descent Directions}, we present an algorithm to efficiently compute such descent directions (under the assumption that at least one subderivative can be computed at every point) and prove its feasibility. Using this algorithm, we introduce a descent method for locally Lipschitz continuous MOPs in general Hilbert spaces (Algorithm~\ref{algo:nonsmooth_descent_method}) in Subsection~\ref{subsec:descent_method_nonsmooth_MOPs}. We prove that this method generates sequences of iterates with Pareto critical accumulation points in Theorem~\ref{thm:convergence_proof_nonsmooth_descent_method}. In Section~\ref{sec:Application}, we demonstrate and analyze the behavior of our method in application to a multiobjective obstacle problem on a two-dimensional domain. Concluding, we summarize our results in Section~\ref{sec:conclusion_and_outlook}.

\section{Theoretical background}
\label{sec:Theoretical Background}

In this section we present the fundamentals to state a necessary optimality condition for nonsmooth multiobjective optimization problems in infinite-dimensional Hilbert spaces using the $\varepsilon$-subdifferential. In scalar nonsmooth optimization the Clarke subdifferential \cite{C1983} and the $\varepsilon$-subdifferential \cite{G1977,K2010} are well-known tools to state optimality conditions and formulate optimization methods. In \cite{MKW2014} optimality conditions for nonsmooth multiobjective optimization are derived in the finite-dimensional case using the Clarke subdifferential. In \cite{G2021} the $\varepsilon$-subdifferential is used for nonsmooth multiobjective optimization problems for finite-dimensional problems. In \cite{M2018,M2019} properties of so-called set-valued gradients are studied in potentially infinite-dimensional Hilbert spaces. These set-valued gradients can be seen as a generalization of the $\varepsilon$-subgradient. The works \cite{M2018, M2019} only focus on the case of scalar optimization.

After introducing the above concepts in this section, we combine them in Section~\ref{sec:eps-subgradient for MOPs} to state a necessary optimality condition for nonsmooth multiobjective optimization problems in infinite-dimensional Hilbert spaces utilizing $\varepsilon$-subdifferentials. We use this optimality condition to prove convergence of an efficient descent method for nonsmooth multiobjective optimization problems in infinite-dimensional Hilbert spaces in Section \ref{sec:Descent_method}.

\subsection{Notations}
\label{subsec:Notatios}

The inner product on a real Hilbert space $\H$ is denoted by $\langle \cdot\,, \cdot \rangle$ with induced norm $\lVert \cdot \rVert \coloneqq \sqrt{\langle\cdot,\cdot\rangle}$.
The topological dual space to $\H$ is denoted by $\H\dual$ and unless otherwise stated, we consider the corresponding dual scalar product $\langle \cdot\,, \cdot \rangle_*\coloneqq \langle R^{-1} (\cdot), R^{-1} (\cdot) \rangle$ and its induced natural norm on the dual space $\| \cdot \|_*$, where $R : \H \to \H\dual$ is the Riesz representation operator.
The symbols $B_{\varepsilon}(x) \coloneqq \{ y \in \H : \| x - y \| < \varepsilon \}$ and $\overline{B_{\varepsilon}(x)} \coloneqq \{ y \in \H : \| x - y \| \le \varepsilon \}$ denote the open and closed $\varepsilon$-balls in $\H$ centered at $x$, respectively. 

For an arbitrary subset $A \subseteq \H$, the symbol $\conv (A)$ is the convex hull of $A$ and $\closure{\conv}(A)$ its closure. 
We denote the $k$-dimensional positive unit simplex by $\Delta^k \coloneqq \{ \lambda \in \R^k : \sum_{i=1}^k \lambda_i = 1, \, \lambda_i \ge 0 \, \text{ for all } \, i = 1,\dots, k\}$.

\subsection{Nonsmooth multiobjective optimization}
\label{sec:introduction}

Let us consider the following nonsmooth multiobjective optimization problem
\begin{align} \tag{\textbf{MOP}} \label{eq:MOP}
	\min_{x \in \H} f(x) = \min_{x \in \H}
	\begin{pmatrix}
		f_1(x) \\
		\vdots \\
		f_k(x)
	\end{pmatrix},
\end{align}
where $f :\H \rightarrow \R^k$ is the \emph{objective vector} with the \emph{objective functions} $f_i :\H \rightarrow \R$ for $i = 1,\dots,k$.

Recall that a function $f:\H\to\R$ is called \emph{locally Lipschitz} near $x\in\H$, if there exist $\varepsilon>0$ and a constant $L=L(x,\varepsilon) > 0$ with
\begin{align*}
	|f(y)- f(z)|\leq L\, \|y-z\|\quad\text{for all }y, z\in B_{\varepsilon}(x).
\end{align*}
We say that $f$ is \emph{locally Lipschitz of rank $L$} if we want to point out the specific Lipschitz constant.

Since \eqref{eq:MOP} is an optimization problem with a vector-valued objective function, the classical concept of optimality from the scalar case cannot be conveyed directly. Instead, we are looking for the \emph{Pareto set}, which is defined in the following way:

\begin{definition}[{\cite[pp. 10-20]{Miettinen1998}}]
    Consider the optimization problem \eqref{eq:MOP}.
    \begin{enumerate}[i)]
        \item [\rm a)] A point $x^* \in \H$ is \emph{Pareto optimal} if there does not exist another point $x \in \H$ such that $f_i(x) \le f_i(x^*)$ for all $i = 1,\dots,k,$ and $f_j(x) < f_j(x^*)$ for at least one index $j$. The set of all Pareto optimal points is the \emph{Pareto set}, which we denote by $P$. 
        \item [\rm b)] A point $x^* \in \H$ is \emph{weakly Pareto optimal} if there does not exist another point $x \in \H$ such that $f_i(x) < f_i(x^*)$ for all $i = 1,\dots,k$.
    \end{enumerate}
    If there exists a $\delta>0$ such that the respective conditions in the definitions above hold only for all $x\in B_{\delta}(x^*)$, then $x^*$ is called \emph{locally} (weakly) Pareto optimal.
\end{definition}
In practice, to check if a given point is Pareto optimal, we need optimality conditions. In the smooth case, there are the well-known \emph{Karush-Kuhn-Tucker (KKT) conditions} (cf.~\cite{Miettinen1998}, for instance), which are based on the gradients of the objective functions. If the objective functions are merely locally Lipschitz, the KKT conditions can be generalized using the concept of \emph{subdifferentials}. In the following, we recall the required definitions and results from nonsmooth analysis. For a more detailed introduction, we refer to \cite[Chapter~2]{C1983}.	

\subsection{Generalized gradients}

\begin{definition}\label{def:gen_dir_deriv}
    For $f:\H\to \R$ locally Lipschitz define the \emph{generalized directional derivative} at $x$ in a direction $v\in\H$ as
    \begin{align}
        f^{\circ }(x;v) \coloneqq \limsup_{y \to x, t \searrow 0} \frac{f(y+t v)-f(y)}{t}.
    \end{align}
\end{definition}
    
In the following we refer to Propositions~2.1.1, 2.1.2 and 2.1.5 in \cite{C1983} which state the most important facts on the generalized directional derivative for our analysis.
    \clearpage
\begin{proposition}\label{prop:gen_dir_deriv}
    Let $f:\H\to \R$ be locally Lipschitz of rank $L$ near $x\in\H$. Then:
    \begin{itemize}
        \item[\em a)] The function $v\mapsto f^\circ(x;v)$ is finite, positively homogeneous, and subadditive on $\H$ (i.e., $f^\circ(x;tv)=tf^\circ(x;v)$ and $f^\circ(x;v+w)\le f^\circ(x;v)+f^\circ(x,w)$ for every $t>0$ and $v,w\in\H$), and satisfies
        \begin{align*}
            \lvert f^{\circ }(x;v) \rvert \le L\,\lVert v \rVert.
        \end{align*}
        \item[\em b)] $f^{\circ }(x;v)$ is upper semicontinuous as a function of $(x,v)$ and, as a function of $v$ alone, is Lipschitz of rank $L$ on $\H$.
        \item[\em c)] $f^{\circ }(x;-v) = (-f)^{\circ }(x;v)$.
    \end{itemize}
\end{proposition}
    
Using the generalized directional derivative we are able to define the so-called \emph{(Clarke) subdifferential}.

\begin{definition} \label{def:Clarke_subdiff}
    For $f:\H\to \R$ locally Lipschitz define the \emph{(Clarke) subdifferential} at $x$ as
	\begin{align*}
		\partial f(x) \coloneqq\big\{\xi \in \H\dual : f^{\circ }(x;v) \ge \xi(v) \text{ for all } v \in \H\big\}.
	\end{align*}	
	A functional $\xi$ in the set $\partial f(x)$ is called a \emph{subderivative} of $f$ in $x$.
\end{definition}
	
If an objective function is continuously (Fréchet-)differentiable, the Clarke subdifferential is a singleton containing only the derivative.
	
\begin{proposition}
    \label{prop:partial_weak_compacrt}
	Let $f:\H\to \R$ be locally Lipschitz of rank $L$ near $x\in\H$. Then:
	\begin{itemize}
	    \item[\em a)] $\partial f(x)$ is a nonempty, convex, weakly compact subset of $\H^*$ and $\lVert \xi \rVert_* \le L$ for every $\xi$ in $\partial f(x)$.
	    \item[\em b)] For every $v$ in $\H$, one has
        \begin{align}
	        f^{\circ }(x;v) = \max\big\{ \xi(v) : \xi \in \partial f(x)\big\}.
	    \end{align}
     \end{itemize}
\end{proposition}
	
\begin{proposition}
	\label{prop:partial_weak_closed}
	Let $f:\H\to \R$ be locally Lipschitz near $x$. Then:
	\begin{itemize}
	    \item[\em a)] We have $\xi \in \partial f(x)$ if and only if $f^{\circ }(x;v) \ge \xi(v)$ for all $v \in \H$.
	    \item[\em b)] Let $(x_i)_i$ and $(\xi_i)_i$ be sequences in $\H$ and $\H^*$, respectively, with $\xi_i \in \partial f(x_i)$. Suppose that $(x_i)_i$ converges to $x$ and that $\xi$ is a weak accumulation point of $(\xi_i)_i$. Then $\xi \in \partial f(x)$.
	    \item[\em c)] $\partial f(x) = \bigcap_{\varepsilon > 0} \bigcup_{y\in \overline{B_{\varepsilon}(x)}} \partial f(y)$.
	\end{itemize}
\end{proposition}
	
\begin{remark}
    In the next subsection we prove modified versions of parts b) and c) of Proposition \ref{prop:partial_weak_closed}. These extensions are required to show our main convergence result.\hfill$\Diamond$
\end{remark}
	
Recall that the Clarke subdifferential in infinite dimensions satisfies the well-known \emph{mean value theorem} (cf., e.g., \cite[Theorem 2.3.7]{C1983}).
	
\begin{theorem}\label{thm:mean_value_thm}
	Let $x$ and $y$ be points in $\H$, and suppose that $f:\H\to \R$ is Lipschitz on an open set containing the line segment $[x,y]$. Then, there exists a point $z$ on the open line segment $(x,y)$ such that
	\begin{align*}
	    f(y)-f(x)\in\partial f(z)(y-x).
	\end{align*}
\end{theorem}
	
Note that, if $f$ is locally Lipschitz continuous on $\H$, then any line segment $[x,y]$ has a neighborhood on which $f$ is globally Lipschitz since $[x,y]$ is compact in $\H$. Using the subdifferential, we can state a necessary optimality condition for locally Lipschitz MOPs.
	
\begin{theorem} \label{thm:KKT}
	Let $f=(f_1,\ldots,f_k):\H\to \R^k$ and $x \in \H$ be a (locally weak) Pareto optimum. Then:
	\begin{align} \label{eq:KKT}
		0 \in \conv \left( \bigcup_{i=1}^k \partial f_i(x) \right).
	\end{align}
	If a vector $x$ satisfies \eqref{eq:KKT} we call it \emph{Pareto critical}.
\end{theorem}
	
\begin{proof}
	We can argue as in \cite[Theorem~12]{MEK2014}. Notice that in \cite{MEK2014} the finite-dimensional case ($\H=\R^n$) is considered. However, the proof can be applied in the infinite-dimensional setting as well without any adjustments. The arguments only rely on properties of the generalized directional derivative and the Clarke subdifferential that we stated in Propositions~\ref{prop:gen_dir_deriv}, \ref{prop:partial_weak_closed} and Theorem~\ref{thm:mean_value_thm} above. 
\end{proof}
	
\begin{remark}
	In the smooth case, \eqref{eq:KKT} reduces to the well-known classical multiobjective KKT conditions. However, in contrast to the smooth case, the optimality con\-di\-tion \eqref{eq:KKT} is numerically challenging to work with, as subdifferentials are dif\-fi\-cult to compute. Therefore, in numerical methods, \eqref{eq:KKT} is only used implicitly.\hfill$\Diamond$
\end{remark}
	
In the following section, we will describe a new way to compute descent directions for nonsmooth MOPs by systematically computing an approximation of $\conv (\cup_{i = 1}^k \partial f_i(x))$ that can be used to obtain a `sufficiently good' descent direction. In addition we use these notions to define `approximate' Pareto critical points in Definition \ref{def:eps_subdiff_mop} which are more stable than actual Pareto critical points from a numerical point of view.
	
\subsection{$\bm\varepsilon$-subdifferentials}

In finite dimensions, $\partial f(x)$ is the convex hull of the limits of the derivatives of $f$ in all sequences (where the derivatives are defined) near $x$ that converge to $x$. Thus, if we evaluate $f'$ in a number of points close to $x$ (where it is defined) and take the convex hull, we expect the resulting set to be an approximation of $\partial f(x)$. To formalize this, we introduce the following definition (cf. \cite{G1977,K2010}).
	
\begin{definition} \label{def:eps_subdiff}
	Let $f:\H\to \R$, $\varepsilon \geq 0$ and $x \in \H$. Then
    \begin{align*}
		\partial_\varepsilon f(x) := \closure{\conv} \Bigg( \bigcup_{y \in \overline{B_\varepsilon(x)}} \partial f(y) \Bigg)
	\end{align*}
	is the closure of the \emph{(Goldstein)} $\varepsilon$-\emph{subdifferential of} $f$ \emph{in} $x$. We call $\xi \in \partial_\varepsilon f(x)$ an $\varepsilon$-\emph{subderivative}.
\end{definition}
	
Note that $\partial_0 f(x) = \partial f(x)$ and $\partial f(x) \subseteq \partial_\varepsilon f(x)$ for all $\varepsilon > 0$. 
	
\begin{proposition} \label{prop:Goldstein}
	Let $x\in\H$ be given arbitrarily and $f:\H\to \R$ be globally Lipschitz continuous on the ball $B_{\bar\varepsilon}(x)$ for some $\bar\varepsilon > 0$. Moreover, suppose that $\varepsilon\in[0,\bar\varepsilon)$. Then $\partial_\varepsilon f(x)$ is nonempty, convex and weakly compact.
\end{proposition}
	
\begin{proof}
	For $\partial_\varepsilon f(x)$, the claim was shown in \cite[Proposition~2.3]{M2018}. To apply the proof we need a neighbourhood of $\overline{B_{\varepsilon}(x)}$, where $f$ is globally Lipschitz continuous. For that reason we introduce the open ball $B_{\bar\varepsilon}(x)\supsetneq \overline{B_{\varepsilon}(x)}$ in the formulation of this proposition. 
\end{proof}
	
In the following, we present a theorem that is a stronger version of parts b) and c) of Proposition \ref{prop:partial_weak_closed}.  This result relates the $\varepsilon$-subdifferential to the Clarke subdifferential. Before we state the theorem we prove a preparatory lemma.
	
\begin{lemma}\label{lem:eps_subdiff_sep}
	Let $f:\H\to\R$ be locally Lipschitz near $x \in \H$, $v \in \H\setminus \{0\}$ and $\alpha \in \R$. If 
	\begin{align}
	    \label{eq:subdiff_sep}
	    \alpha > \xi(v) \quad \text{  for all  } \quad \xi \in \partial f(x),
	\end{align}
	then there exists an $\overline{\varepsilon} > 0$, such that for all $0 \le \varepsilon \le \overline{\varepsilon}$
	\begin{align}
	    \label{eq:eps_subdiff_sep}
	    \alpha > \xi(v)\quad \text{  for all  } \quad \xi \in \partial_{\varepsilon} f(x).
	\end{align}
\end{lemma}
	
\begin{remark}
    Lemma~\ref{lem:eps_subdiff_sep} states that the $\varepsilon$-subdifferential contracts in a well-behaved manner to the Clarke subdifferential as $\varepsilon\searrow0$. In view of Proposition~\ref{prop:partial_weak_closed}-c) this lemma states that we do not have to take the full set-valued limit to contract the $\varepsilon$-subdifferential to one side of the hyperplane defined by $v$ and $\alpha$.\hfill$\Diamond$
\end{remark}
	
\begin{proof}[Proof of Lemma~{\em\ref{lem:eps_subdiff_sep}}]
	We do not show \eqref{eq:eps_subdiff_sep} directly but first conclude that the separation holds in the weaker form of
	\begin{align}
	    \label{eq:eps_ball_sep}
	    \alpha > \xi(v) \quad \text{  for all  } \quad \xi \in \bigcup_{y\in \overline{B_{\varepsilon}(x)}} \partial f(y)\subset\partial_\varepsilon f(x),
	\end{align}
	which is a consequence of Proposition \ref{prop:partial_weak_closed} as we prove in the following.\hfill\\
	Let $v \in \H\setminus \{0\}$ and $\alpha \in \R$. Assume for all $\bar\varepsilon > 0$ there exists an $\varepsilon\in(0,\bar\varepsilon]$ and $\xi \in \bigcup_{y\in \overline{B_{\varepsilon}(x)}} \partial f(x)$ with $\xi(v) \ge \alpha$. Then, there exist a sequence $(\varepsilon_i)_i$ of positive real numbers and sequences $(y_i)_i$ and $(\xi_i)_i$ of elements in $\H$ and $\H^*$, respectively, such that $(\varepsilon_i)_i$ converges to zero, $\xi_i \in \partial f(y_i)$, $\xi_i(v) \ge \alpha$ for all $i\ge 1$ and $\lVert y_i- x \rVert < \varepsilon_i$ converges to zero. Since $f$ is locally Lipschitz continuous in $x$, there exists an $I\ge 1$ such that for all $i \ge I$ the mapping $f$ is locally Lipschitz continuous of rank $L$ in $y_i$. Then, Proposition \ref{prop:partial_weak_compacrt}-a) states that for all $i \ge I$ the elements of the sequence $(\xi_i)_{i}$ are contained in the weakly compact set $\overline{B_{L}(0)}$. Therefore, the sequence $(\xi_i)_{i}$ has a weak sequential accumulation point $\xi^*$. By Proposition \ref{prop:partial_weak_closed}-b) the point $\xi^*$ is an element of $\partial f(x)$. Since $\xi_i(v) \ge \alpha$ for all $i\ge 1$ we get by the weak convergence of a subsequence of $(\xi_i)_i$ to $\xi^*$, that $\xi^*(v) \ge \alpha$ which is a contradiction to \eqref{eq:subdiff_sep}. Therefore, \eqref{eq:eps_ball_sep} holds.\hfill\\
	The remainder of the proof follows by the definition of the $\varepsilon$-subdifferential (cf. Definition~\ref{def:eps_subdiff}). If a set lies on one side of a hyperplane, then also its convex hull lies on that side and also its closure. 
\end{proof}
	
\begin{theorem}
	\label{thm:subdiff_eps_subdiff}
	    Let $f:\H\to\R$ be locally Lipschitz near $x \in \H$. Then, the following statements hold:
	    \begin{itemize}
	        \item[\em a)] Let $(x_i)_i$ be a sequence in $\H$ converging to $x$ and $(\varepsilon_i)_i$ a sequence in $\R_{> 0}$ tending to $0$. Suppose that the sequence $(\xi_i)_i$ satisfies $\xi_i\in\partial_{\varepsilon_i} f(x_i)$. Let $\xi$ be a weak sequential accumulation point of $(\xi_i)_i$. Then $\xi \in \partial f(x)$.
	        \item[\em b)] $\partial f(x) = \bigcap\limits_{\varepsilon > 0} \partial_{\varepsilon} f(x)$.
	    \end{itemize}
\end{theorem}
	
\begin{proof}
	\begin{itemize}
        \item [a)] Since $\xi_i \in \partial_{\varepsilon_i} f(x_i)$ it follows that $\xi_i \in \partial_{\kappa_i} f(x)$, with $\kappa_i = \varepsilon_i + \lVert x_i - x \rVert$. Assume $\xi \notin \partial f(x)$. Then, since $\partial f(x)$ is convex and weakly compact, it is closed and the strict separation theorem states that there exists $v\in \H\setminus \{0\}$ and $\alpha \in \R$ satisfying
	    \begin{align*}
	       \xi(v) > \alpha > \eta(v) \quad \text{for all} \quad \eta \in \partial f(x).
	    \end{align*}
	    Since $\kappa_i$ converges to $0$, Lemma \ref{lem:eps_subdiff_sep} states that there exists an $I \ge 1$ such that
	    \begin{align*}
	       \xi(v) > \alpha > \eta(v) \quad \text{for all} \quad \eta \in \partial_{\kappa_i} f(x) , \quad i \ge I,
	    \end{align*}
        and hence
        \begin{align*}
	       \xi(v) > \alpha > \xi_i(v) \quad \text{for all}\quad i \ge I.
	    \end{align*}
        This is a contradiction to the fact that $\xi$ is a weak sequential accumulation point of $(\xi_i)_i$.
	    \item [b)] From Proposition~\ref{prop:partial_weak_closed}-c) we immediately get the inclusion
	    \begin{align*}
	        \partial f (x) & = \bigcap_{\varepsilon > 0} \bigcup_{y \in \overline{B_{\varepsilon}(x)}} \partial f(y)\subseteq \bigcap_{\varepsilon > 0} \overline{\conv}\bigg(\bigcup_{y \in \overline{B_{\varepsilon}(x)}} \partial f(y)\bigg)= \bigcap_{\varepsilon > 0} \partial_{\varepsilon}f(x).
	    \end{align*}
	    The other inclusion is a consequence of Lemma \ref{lem:eps_subdiff_sep} and we prove it analogously to part a): Assume that $\xi \in\cap_{\varepsilon > 0} \partial_{\varepsilon} f(x)$, but $\xi\notin \partial f(x)$. Then, since $\partial f(x)$ is convex and weakly compact and therefore closed, the strict separation theorem states that there exist $v\in \H\setminus \{0\}$ and $\alpha \in \R$ such that
        \begin{align*}
	        \xi(v) > \alpha > \eta(v) \quad \text{for all} \quad \eta \in \partial f(x).
	    \end{align*}
        Lemma \ref{lem:eps_subdiff_sep} states that there exists an $\varepsilon > 0$ such that
        \begin{align*}
    	    \xi(v) > \alpha > \eta(v) \quad \text{  for all  } \quad \eta \in \partial_{\varepsilon} f(x)
	    \end{align*}
        and hence $\xi \notin \partial_{\varepsilon} f(x)$. Therefore, it follows that $\xi\notin \bigcap_{\varepsilon > 0} \partial_{\varepsilon} f(x)$, which is a contradiction. In total, we derive $\cap_{\varepsilon > 0} \partial_{\varepsilon} f(x)\subseteq\partial f(x)$ which completes the proof.
    \end{itemize}
\end{proof}

\section{The $\bm\varepsilon$-subdifferential for MOPs}
\label{sec:eps-subgradient for MOPs}

In this section we extend the Goldstein $\varepsilon$-subdifferential to the multiobjective setting. We define a multiobjective $\varepsilon$-subdifferential and investigate its main properties.
    
\begin{definition}\label{def:eps_subdiff_mop}
	Let $f=(f_1,\ldots,f_k):\H\to\R^k$, $\varepsilon \ge 0$ and $x\in \H$. Then
	\begin{align*}
	    F_{\varepsilon}(x) = \overline{\conv}\left( \bigcup_{i = 1}^k \partial_{\varepsilon}f_i(x) \right)
	\end{align*}
	generalizes the $\varepsilon$-subdifferential to the multiobjective setting.
\end{definition}

We use the multiobjective $\varepsilon$-subdifferential to give an approximate notion of criticality with the following definition.

\begin{definition}\label{def:eps_delat_crit}
    We say that $x \in \H$ is $(\varepsilon, \delta)$\emph{-critical} for constants $\varepsilon \geq 0$ and $\delta \geq 0$, if there exists a $\xi \in F_{\varepsilon}(x)$ with $\lVert \xi \rVert_* \le \delta$.
\end{definition}
		
\begin{lemma}\label{Lemma_WeaklyCompact}
	The convex hull of a finite union of convex, weakly compact sets is weakly compact.
\end{lemma}
	
\begin{proof}
	Although the proof utilizes standard arguments, we state it here for the sake of completeness. Let $A^i \subseteq\H$ be nonempty, convex and weakly compact for all $i \in \{1,\dots,k\}$ and set $A = \conv( \cup_{i=1}^k A^i)$. Let $(x_m)_m$ be an arbitrary sequence in $A$. Since the sets $A^i$ are convex, we can write 
	\begin{align*}
	    x_m=\sum_{i=1}^k\lambda_m^i x_m^i \quad \text{for all }m,
	\end{align*}
	with $\lambda_m = (\lambda_m^1, \dots, \lambda_m^k)^{\top} \in \Delta^k$ and $x_m^i \in A^i$. Since $\Delta^k$ is compact and the sets $A^i$ are weakly sequentially compact, there exists a subsequence $(m_l)_l$ such that $\lambda_{m_l}$ converges to $\lambda_* \in \Delta^k$ and that the subsequences $(x_{m_l}^i)_l$ converge weakly to $x^i_* \in A^i$ for $i=1,\ldots,k$. Then $(x_{m_l})_l$ converges weakly to $x_*=\sum_{i=1}^k \lambda_*^i x_*^i \in \conv( \cup_{i = 1}^k A_i)$, which completes the proof.  
\end{proof}
 
Now, we formulate the following result analogously to Propositon~\ref{prop:Goldstein}.
	
\begin{proposition}\label{prop:F_eps}
	For $i = 1,\dots, k$ let $f_i:\H\to\R$ be globally Lipschitz on $B_{\bar\varepsilon}(x)$ for some $x\in\H$ and $\bar\varepsilon > 0$ and let $\varepsilon\in[0,\bar\varepsilon)$. Then $F_{\varepsilon}(x)$ is nonempty, convex and weakly compact. Furthermore,
    \begin{align}\label{eq:F_eps_wo_wclosure}
	    F_{\varepsilon}(x) = \conv\left( \bigcup_{i = 1}^k \partial_{\varepsilon}f_i(x) \right),
	\end{align}
    i.e., the closure in Definition~{\em\ref{def:eps_subdiff_mop}} is superfluous in this case.
\end{proposition}
	
\begin{proof}
	The proof follows from Proposition~\ref{prop:Goldstein} and Lemma~\ref{Lemma_WeaklyCompact}.
\end{proof}
	
The following theorem extends Theorem~\ref{thm:subdiff_eps_subdiff} to the multiobjective setting.
	
\begin{theorem}
    \label{thm:subdiff_eps_subdiff_multi}
	For $i = 1,\dots, k$ let $f_i:\H\to\R$ be locally Lipschitz near $x \in \H$. Let $(\varepsilon_j)_j$ be a sequence of positive numbers that converges to 0. Let $(x_j)_j$ and $(\xi_j)_j$ be sequences in $\H$ and $\H^*$, respectively, and assume that $(x_j)_j$ converges to $x$ and that $(\xi_j)_j$ tends $\H^*$-weakly to $\xi$. Further assume that $\xi_j \in F_{\varepsilon_j}(x_j)$ for all $j \ge 1$. Then,
	\begin{align}
	    \xi \in F_{0}(x) = \conv\left( \bigcup_{i=1}^k \partial f_i(x)\right).
	\end{align}
\end{theorem}
	
\begin{proof}
    Since the functions $f_i$ are locally Lipschitz continuous for $i = 1,\dots, k$, there exists $\varepsilon > 0$ such that all $f_i$ are Lipschitz continuous of rank $L$ on $B_\varepsilon(x)$. Similar to the proof of Theorem \ref{thm:subdiff_eps_subdiff} we define $\kappa_j = \varepsilon_j + \lVert x_j - x \rVert$ and fix $J \ge 1$ such that for all $j \ge J$ it holds that $\kappa_j \le \varepsilon$. From $\partial_{\varepsilon_j}f(x_j) \subseteq \partial_{\kappa_j}f(x)$ it follows that $F_{\varepsilon_j}(x_j) \subseteq F_{\kappa_j}(x)$. Proposition \ref{prop:F_eps} implies that $F_{\kappa_j}(x)$ is nonempty, convex and weakly compact and
    \begin{align*}
        F_{\kappa_j}(x) = \conv \left(\bigcup_{i=1}^k \partial_{\kappa_j}f_i(x)\right).
    \end{align*}
    The remainder of the proof can be seen as a combination of the proofs of Theorem \ref{thm:subdiff_eps_subdiff} and Proposition \ref{prop:F_eps}. Since $\xi_j$ is an element of $F_{\kappa_j}(x)$ for all $j \ge J$ it can be written as
    \begin{align*}
        \xi_j = \sum_{i=1}^k \lambda_j^i \xi_j^i,
    \end{align*}
    with $\lambda_j = (\lambda_j^1,\dots \lambda_j^k) \in \Delta^k$ and $\xi_j^i \in \partial_{\kappa_j} f_i(x)$. Since $\kappa_j \le \varepsilon$ it follows that $\xi_j^i$ is contained in the weakly compact set $\overline{B_L(0)}$. Hence, there exists a subsequence $(j_l)_l$ such that
    \begin{align*}
        \lambda_{j_l} \to \lambda_* \in \Delta^k\, \text{ and }\, \xi_{j_l}^i \rightharpoonup \xi_*^i \in \overline{B_L(0)}\text{ for  all }i=1,\ldots,k.
    \end{align*}
    From Theorem \ref{thm:subdiff_eps_subdiff} it follows that $\xi_*^i \in \partial f_i(x)$. Then $\xi_j = \sum_{i=1}^k \lambda_j^i \xi_j^i$ converges weakly to $\xi_*=\sum_{i=1}^k \lambda_*^i \xi_*^i \in \conv( \cup_{i=1}^k \partial f_i(x))$. Since the weak limit is unique and $\xi_j$ converges weakly to $\xi$ the proof is complete. 
\end{proof}
	
The next corollary follows directly from Theorem \ref{thm:subdiff_eps_subdiff_multi} and gives a sufficient condition for a point to be Pareto critical.
	
\begin{corollary}\label{cor:subdiff_eps_subdiff_multi}
	For $i = 1,\dots, k$ let $f_i:\H\to\R$ be locally Lipschitz in $x\in\H$. Assume that 
    \begin{align*}
	    0 \in F_\varepsilon(x) \quad \text{ for all } \varepsilon > 0.
	\end{align*}
    Then x is Pareto critical, i.e.,
    \begin{align*}
	    0 \in \conv \left( \bigcup_{i=1}^k \partial f_i(x) \right).
	\end{align*}
\end{corollary}

After describing the optimality conditions for MOPs, we now move towards the algorithms from \cite{G2021} that we extend to the infinite-dimensional setting.

\section{Descent method for nonsmooth MOPs}	
\label{sec:Descent_method}

In this section, we present a line-search based \emph{common-descent method}, meaning that, starting from a point $x_1 \in \H$, we generate a sequence $(x_j)_j$ in $\H$ in which each point is an improvement over the previous point with respect to all objective functions, i.e.,
\begin{align*}
	f_i(x_{j+1}) < f_i(x_j) \quad \text{for all }j \ge 1\text{ and }i = 1,\dots,k,
\end{align*}
and where $x_{j+1} = x_j + t_j v_j$ for a search direction $v_j\coloneqq R^{-1}(\xi_j)$ generated from a dual element $\xi_j \in \H^*$ and corresponding step lengths $t_j \in \R_{>0}$. The critical computation of the search direction generalizes the method from \cite{G2021} to the infinite-dimensional setting.

The foundation of our approach is the following result from convex analysis.

\begin{theorem} \label{thm:steepest_descent_direction}
	Let $\WW \subseteq\H\dual$ be convex and closed. Then,
	\begin{align} \label{eq:min_norm_problem}
		\bar \xi\coloneqq\argmin_{\xi \in -\WW}{\|\xi\|}_*^2
	\end{align}		
	is well-defined and unique. Further, it holds that either $\bar{\xi} \neq 0$ and
	\begin{align}
        \label{eq:steepest_descent_ineq}
		{\langle\bar \xi,\xi\rangle}_*\leq-{\|\bar \xi\|}_*^2<0\quad\text{for all } \xi \in \WW,
	\end{align}	
	or $\bar{\xi} = 0$ and there is no $\tilde \xi \in \H$ with $\langle \tilde \xi, \xi \rangle_* < 0$ for all $\xi \in \WW$.
\end{theorem}
	
\begin{proof}
	This theorem is stated in \cite[Theorem~3.14]{B2011}. 
\end{proof}

When considering $\WW = F_\varepsilon(x)$ (which is convex and closed by definition), then this immediately yields the following corollary.
    \clearpage
\begin{corollary}\label{cor:pareto_opt}
	Let $\varepsilon \geq 0$.
	\begin{itemize}
		\item[\em a)] If $x$ is locally weakly Pareto optimal, then
		\begin{align} \label{eq:eps_critical}
			0 \in F_\varepsilon(x).
		\end{align}
		\item[\em b)] Let $x \in \H$ and
		\begin{align}
            \label{eq:eps_min_norm_problem}
			\bar \xi\coloneqq\argmin_{\xi \in -F_\varepsilon(x)}{\|\xi\|}_*^2.
		\end{align}		
		Then either $\bar{\xi} \neq 0$ and
		\begin{align} \label{eq:eps_descent_ineq}
			{\langle\bar \xi,\xi\rangle}_*\leq-{\|\bar \xi\|}_*^2 < 0\quad\text{for all }\xi \in F_\varepsilon(x),
		\end{align}	
		or $\bar \xi=0$ and there is no $\tilde\xi \in \H$ with $\langle \tilde\xi, \xi \rangle_* < 0$  for all $\xi \in F_\varepsilon(x)$.
	\end{itemize}
\end{corollary}	

This means that, when working with the $\varepsilon$-subdifferential instead of the Clarke subdifferential, we still have a necessary optimality condition and a way to compute descent directions, although the optimality conditions are weaker and descent can be expected to be weaker than when using the unrelaxed subdifferential.

For the direction from \eqref{eq:eps_min_norm_problem}, we can find a lower bound for a step size up to which we have guaranteed descent in each objective function $f_i$.

\begin{lemma} \label{lem:step_size_bound}
	For $i = 1,\dots, k$ let $f_i\colon\H\to\R$ be locally Lipschitz continuous in $x\in\H$. Moreover, we assume that $\varepsilon \geq 0$ holds and we define $\bar v\coloneqq R^{-1}(\bar \xi)$ for the solution $\bar{\xi} \in - F_{\varepsilon}(x)$ of \eqref{eq:eps_min_norm_problem}. Then
    \begin{align*}
		f_i(x + t \bar{v}) \leq f_i(x) - t\,\lVert \bar{v} \rVert^2 \quad \text{for all }0 \leq t \leq \frac{\varepsilon}{\lVert \bar{v} \rVert}\text{ and }i \in \{1,\dots,k\}.
	\end{align*}
\end{lemma}
	
\begin{proof}
	The proof of \cite[Lemma~3.2]{G2021} can be adapted to the infinite-dimensional case using the fact that the mean value theorem (Theorem \ref{thm:mean_value_thm}) holds for the Clarke subdifferential also in infinite dimensions and because $\|\bar v\|=\|\bar\xi\|_*$. 
\end{proof}
    
However, solving \eqref{eq:min_norm_problem} generally requires the knowledge of the \emph{entire} $\varepsilon$-sub\-differential, which is impractical. Instead, we will use Theorem \ref{thm:steepest_descent_direction} to compute a finitely generated approximation $\WW$ of $\conv \left( \cup_{i = 1}^k \partial f_i(x) \right)$, where the resulting direction is guaranteed to have sufficient descent.
	
\subsection{Efficient computation of descent directions}
\label{subsec:Efficient Computation of Descent Directions}

In this subsection, we describe how the solution of \eqref{eq:eps_min_norm_problem} can be replaced by a suboptimal one when only a single subderivative is available in every $x \in \H$. Similar to the gradient sampling approach, the idea behind this method is to use instead of $F_\varepsilon(x)$ in \eqref{eq:eps_min_norm_problem} the convex hull of a finite number of $\varepsilon$-subderivatives $\xi_1,\dots,\xi_m$ from $F_\varepsilon(x)$ for $m \ge 1$. Since it is impossible to know a priori how many and which $\varepsilon$-subderivatives are required to obtain a good descent direction, we solve \eqref{eq:eps_min_norm_problem} multiple times in an iterative manner while enriching our approximation until a satisfying direction has been found. To this end, in the following, we will specify how to enrich our current approximation $\conv(\{ \xi_1, \dots, \xi_m \})$ and how to characterize an acceptable descent direction.\newline

Suppose that $\WW = \{\xi_1, \dots, \xi_m\} \subseteq F_\varepsilon(x)$ and define
\begin{align}
    \label{eq:approx_desc_dir}
	\tilde{\xi}:=\argmin_{\xi\in-\conv(\WW)}{\|\xi\|}_*^2.
\end{align}
Let $c \in(0,1)$. Motivated by Lemma \ref{lem:step_size_bound}, we regard $\tilde v \coloneqq R^{-1}(\tilde{\xi})$ as an \emph{acceptable} descent direction, if
\begin{align} \label{eq:accept_direction}
	f_i \left( x + \frac{\varepsilon}{\lVert \tilde{v} \rVert} \tilde v \right) \leq f_i(x) - c \varepsilon \lVert \tilde{v} \rVert \quad \text{for all }i \in \{1,\dots,k\}.
\end{align}
If the set $I \subseteq \{1,\dots,k\}$ for which \eqref{eq:accept_direction} is violated is non-empty then we have to find a new $\varepsilon$-subderivative $\xi' \in F_\varepsilon(x)$ such that $\WW \cup \{ \xi' \}$ yields a better descent direction. Intuitively, \eqref{eq:accept_direction} being violated means that the local behavior of $f_i$, $i \in I$, in $x$ in the direction $\tilde{v}$ is not sufficiently captured in $\WW$. Thus, for each $i \in I$, we expect that there exists some $t'\in(0,\varepsilon/\lVert \tilde{v} \rVert]$ such that $\xi' \in \partial f_i(x + t' \tilde{v})$ improves the approximation of $F_\varepsilon(x)$. This is stated in the following lemma. For a proof, we refer to \cite[Lemma~3.3]{G2021}.
	
\begin{lemma} \label{lem:infeasible_direction}
	Let $c\in(0,1)$, $\WW = \{ \xi_1, \dots, \xi_m \} \subseteq F_\varepsilon(x)$ and $\tilde{v}\coloneqq R^{-1}(\tilde \xi)$ for the solution $\tilde \xi$ of \eqref{eq:approx_desc_dir} and assume $\tilde{v} \neq 0$. If
    \begin{align*}
		f_i \left( x + \frac{\varepsilon}{\lVert \tilde{v} \rVert} \tilde{v} \right) > f_i(x) - c \, \varepsilon \lVert \tilde{v} \rVert \quad\text{for some }i \in \{1,\dots,k\},
	\end{align*}
	then there is some $t'\in(0,\varepsilon/\lVert \tilde{v} \rVert]$ and $\xi' \in \partial f_i(x + t' \tilde{v})$ such that
    \begin{align}
        \label{eq:new_subgrad_condition}
		{\langle\tilde{\xi},\xi'\rangle}_*>-c\,{\|\tilde{\xi}\|}_*^2.
	\end{align}
	In particular, $\xi' \in F_\varepsilon(x) \setminus \conv(\WW)$.
\end{lemma}

Note that Lemma \ref{lem:infeasible_direction} only shows the existence of $t'$ and $\xi'$ without stating a way how to actually compute them. To solve this problem, let $i \in \{1,\dots,k\}$ be the index of an objective function for which \eqref{eq:accept_direction} is not satisfied, define
\begin{align}
    \label{eq:def_h}
    	h_i : \R \rightarrow \R, \quad t \mapsto f_i(x + t \tilde{v}) - f_i(x) + c t\,\lVert \tilde{v} \rVert^2
\end{align}
and consider Algorithm~\ref{algo:new_subgradient}.
\begin{algorithm} 
	\caption{(Computing of a new subderivative)}
	\label{algo:new_subgradient}
	\begin{algorithmic}[1] 
		\REQUIRE Current point $x \in \H$, direction $\tilde{v}=R^{-1}(\tilde \xi)\in\H$, tolerance $\varepsilon>0$, Armijo parameter $c\in(0,1)$.
        \STATE Set $a = 0$, $b=\varepsilon/\lVert \tilde{v}\rVert$ and $t=(a+b)/2$.
        \FOR{$j = 1,2,\dots$}
	        
			\STATE Compute a $\xi' \in \partial f_i(x + t \tilde{v})$.
			\IF{$\langle \tilde{\xi}, \xi' \rangle_* > - c\,\lVert \tilde{\xi} \rVert_*^2$}
                \STATE stop.
            \ENDIF
			\IF{$h_i(b) > h_i(t)$}
                \STATE set $a = t$.
            \ELSE
                \STATE set $b = t$.
            \ENDIF
			\STATE Set $t=(a+b)/2$.
        \ENDFOR
        \RETURN Current $\xi' \in \partial f_i(x + t \tilde{v})$.
    \end{algorithmic} 
\end{algorithm}		
If $f_i$ is continuously differentiable around $x$, then \eqref{eq:new_subgrad_condition} is equivalent to $h_i'(t') > 0$, i.e., $h_i$ being monotonically increasing around $t'$. Thus, the idea of Algorithm \ref{algo:new_subgradient} is to find some $t$ such that $h_i$ is monotonically increasing around $t$, while checking if \eqref{eq:new_subgrad_condition} is satisfied for a subderivative $\xi \in \partial f_i(x + t \tilde{v})$. For a more thorough discussion of the behavior and termination of Algorithm \ref{algo:new_subgradient}, we refer to \cite{G2024,G2021}.

We use this method of finding new subgradients to construct an algorithm that computes descent directions of nonsmooth MOPs, namely Algorithm~\ref{algo:descent_direction}.

\begin{algorithm} 
	\caption{(Computing a descent direction)}
	\label{algo:descent_direction}
	\begin{algorithmic}[1] 
		\REQUIRE Current point $x \in \H$, tolerances $\varepsilon, \delta > 0$, Armijo parameter $c \in(0,1)$.
        \STATE Compute $\xi^i_1 \in \partial_\varepsilon f_i(x)$ for all $i \in \{1,\dots,k\}$. Set $\WW_1 = \{ \xi^1_1, \dots, \xi^k_1 \}$ and $l = 1$.
        \FOR{$l=1,2,\ldots$}
			\STATE Compute $\xi_l = \argmin_{\xi \in -\conv(\WW_l)} \lVert \xi \rVert_*^2$ and set $v_l = R^{-1}(\xi_l)$.
			\IF{$\lVert \xi_l \rVert_* \leq \delta$}
                \RETURN $v_l$.
            \ELSE
			    \STATE Find all objective functions for which there is insufficient descent:
			    \begin{align*}
				    \textstyle I_l = \big\{ j \in \{1,\dots,k\}\,\big|\,f_j \big(x +\varepsilon v_l/\|v_l\| \big) > f_j(x) - c\,  \varepsilon\,\|v_l\|\big\}.
				\end{align*}
				\IF{$I_l= \emptyset$}
                    \STATE stop.
                \ELSE
			        \STATE For each $j \in I_l$ compute $t_j \in(0,\varepsilon/ \lVert v_l \rVert]$ and $\xi^j_l \in \partial f_j(x + t_j v_l)$ with $\langle \xi_l,\xi^j_l \rangle_*> -c\,\|\xi_l\|_*^2$ by applying Algorithm \ref{algo:new_subgradient}.
			        \STATE Set $\WW_{l+1}=\WW_l\cup\{\xi^j_l\,|\,j\in I_l\}$.
                \ENDIF
            \ENDIF
        \ENDFOR
	\end{algorithmic} 
\end{algorithm}

In Theorem~\ref{thm:descent_dir_conv}, we will show that Algorithm~\ref{algo:descent_direction} stops after a finite number of iterations and produces an acceptable descent direction (cf.~\eqref{eq:accept_direction}). In the infinite-dimensional setting, the proof of \cite[Theorem~3.1]{G2021} cannot be applied directly. The proof uses the fact that the closed ball $\overline{B_{\varepsilon}(x)}$ is a compact subset of $\R^n$ to conclude that there exists a common Lipschitz constant $L$ on $B_{\varepsilon}(x)$ for the locally Lipschitz objective functions $f_i$ . This premise does not hold for infinite-dimensional Hilbert spaces. In fact one can construct a function $f$ that is locally Lipschitz on $\H$ but not Lipschitz continuous on $B_2(0)$, as demonstrated in the following example.
	
\begin{example}\label{ex:locally_but_not_globally_lipschitz}
	Let $\H$ be a separable Hilbert space with orthonormal basis $\{e^i\}_{i\in\mathbb N}$. For $i \ge 2$ we define by $\mathcal B_i \coloneqq \overline{B_{1/i}(e_i)}$ a family of closed balls. Obviously, we have $\mathcal B_i \cap \mathcal B_j = \emptyset$ for $i \neq j$, since $\lVert e_i - e_j\lVert = \sqrt{2}>1/i+1/j$. Using the sets $\mathcal B_i$ define the function
	\begin{align*}
	    f:\H\to \R, \quad x\mapsto 
	    \begin{cases*}
            i\,\lVert x-e_i \rVert & if $x \in \mathcal{B}_i$,\\
            1        & otherwise.
        \end{cases*}
	\end{align*}
    The local Lipschitz continuity can be derived from the definition of $f$. In fact, the set $\H \setminus \bigcup_{i\ge 2} \mathcal B_i$ is open and hence for every $x \in \H \setminus \bigcup_{i\ge 2} \mathcal B_i$ there exists a neighborhood of $x$ on which $f$ is constant and therefore Lipschitz continuous. If $x \in \mathcal B_i$ for some $i\ge2$ there exists an open neighboorhood $\mathcal U$ of $x$ such that $\mathcal U \cap \mathcal B_j = \emptyset$ for $j \neq i$. Then for all $y, z \in \mathcal{U}$ it holds that $\lvert f(y) - f(z) \rvert \le i \lVert y - z \rVert$, which can be verified by a simple case seperation considering all the case where $y$ and $z$ belong to $\H \setminus \bigcup_{i\ge 2} \mathcal B_i$ or $\mathcal B_i$.\\
    If $f$ would be Lipschitz continuous on $B_2(0)$ with some Lipschitz constant $L>0$ we arive at a contradiction because then $B_i(0) = \partial f(e_i) \subseteq B_L(0)$ has to hold since $e_i \in B_2(0)$ for all $i \ge 2$.
\end{example}
	
Nevertheless we can show that Algorithm \ref{algo:descent_direction} still converges for an infinite-dimensional Hilbert space. We can recover the main argument of the proof of \cite[Theorem~3.1]{G2021} but need some preparatory results to bypass the fact that we cannot use a common Lipschitz constant for the functions $f_i$ on $B_{\varepsilon}(x)$. To this end, we introduce the two following lemmas.
	
\begin{lemma}\label{lem:convex_update_1}
    Let $C_1 \subseteq C_2 \subseteq \H\dual$ be two convex and closed subsets. Define
    \begin{align*}
    	\xi_1 &\coloneqq \argmin_{\xi \in C_1}{\|\xi\|}_*^2 \quad\text{and}\quad \xi_2 \coloneqq \argmin_{\xi \in C_2}{\|\xi\|}^2_*.
    \end{align*}
    Note that $\xi_1$ and $\xi_2$ are well-defined and unique. Let $\delta \ge 0$  such that $\|\xi_2\|_* \ge \delta$. Then
    \begin{align*}
    	{\|\xi_1-\xi_2\|}_*^2\le{\|\xi_1\|}_*^2-\delta^2.
    \end{align*}
\end{lemma}
	
\begin{proof}
    Simply rewriting the squared norm yields
    \begin{align*}
        {\|\xi_1-\xi_2\|}_*^2={\|\xi_1\|}_*^2-{\|\xi_2\|}_*^2+2\,{\langle\xi_2,\xi_2-\xi_1\rangle}_*.
    \end{align*}
    From $\xi_1 \in C_2$ we infer the projection property $\langle \xi_2, \xi_2 - \xi_1 \rangle_* \le 0$. In addition with the relation $-\lVert \xi_2 \rVert_*^2 \le -\delta^2$ we get the desired result. 
\end{proof}
	
In the proof of the following lemma we directly incorporate Lemma~\ref{lem:convex_update_1}.
	
\begin{lemma}\label{lem:convex_update_2}
	Let $(\xi_l)_l$ be an arbitrary sequence in $\H\dual$. Define $\WW_l \coloneqq \{\xi_1, \dots, \xi_l\}$ for $l \ge 1$. Let the sequence $(\psi_l)_l\subset\H^*$ be given by
	\begin{align*}
	    \psi_l = \argmin_{\psi \in -\conv(\WW_l)}{\|\psi\|}_*^2.
	\end{align*}
	Then $(\psi_l)_l$ converges strongly in $\H\dual$.
\end{lemma}

\begin{proof}
    From the definition of the elements $\psi_l$ we obtain that $( \lVert \psi_l \rVert_* )_l$ is monotonically decreasing. Hence we can conclude that there is a $\delta>0$ such that
	\begin{align*}
	    \lim_{l\to +\infty}{\|\psi_l\|}_*^2=:\delta^2\ge 0.
    \end{align*}
	Using the limit $\delta^2 \ge 0$ and Lemma \ref{lem:convex_update_1}, we will show that $(\psi_l)_l$ is a Cauchy sequence in $\H^*$. Let $l, m \ge 1$ and consider $\|\psi_l-\psi_{l+m}\|_*$. Choosing $C_1=-\conv(\WW_l)$, $C_2 =-\conv(\WW_{l+m})$, $\xi_1=\psi_l$ and $\xi_2=\psi_{l+m}$ with $\|\psi_{l+m}\|_*\ge\delta$ we infer from Lemma~\ref{lem:convex_update_1} that
	\begin{align*}
	    {\|\psi_l-\psi_{l+m}\|}_*^2\le{\|\psi_l\|}_*^2-\delta^2.
    \end{align*}
	Since $\lim_{l \to +\infty} \lVert \psi_l \rVert_*^2 = \delta^2$ it follows that $(\psi_l)_l$ is a Cauchy sequence in $\H\dual$. Consequently, the sequence $(\psi_l)_l$ converges. 
\end{proof}
	
Using Lemmas \ref{lem:convex_update_1} and \ref{lem:convex_update_2} we can adapt the proof of \cite[Theorem~3.1]{G2021} to show that Algorithm \ref{algo:descent_direction} terminates in the Hilbert space setting.
	
\begin{theorem}
    \label{thm:descent_dir_conv}
	For $i = 1,\dots, k$, let $f_i\colon\H\to\R$ be locally Lipschitz continuous. Then, Algorithm~{\em\ref{algo:descent_direction}} terminates so that the sequence $(v_l)_l$ is finite. If $\tilde{v}$ is the last element of $(v_l)_l$ and $\tilde \xi=R(\tilde v)$, then either $\|\tilde \xi\|_* \leq \delta$ or $\tilde v$ is an acceptable descent direction, i.e.,
	\begin{align*}
		f_i \left( x + \frac{\varepsilon}{\lVert \tilde{v} \rVert} \tilde{v} \right) \leq f_i(x) - c \varepsilon \lVert \tilde{v} \rVert \quad \text{for all }i = 1,\dots,k.
	\end{align*}
\end{theorem}
	
\begin{proof}
    Assume Algorithm \ref{algo:descent_direction} does not terminate, i.e., the sequences $(\xi_l)_l$ and $(v_l)_l = \left(R^{-1}(\xi_l)\right)_l$ are infinite sequences. 
    Independently from Steps 7 and 11, Lemma \ref{lem:convex_update_2} guarantees that $(\xi_l)_l$ converges to an element $\tilde{\xi}$ in $\H\dual$, and, accordingly, $(v_l)_l$ converges to $\tilde v = R^{-1}(\tilde \xi)$. Hence, the scalars $t^j_l \in (0,\varepsilon/\lVert v_l \rVert]$ chosen in Step 11 are bounded for all $l \ge 1$ and $j \in I_l$. Using this, we choose a subsequence $(l_m)_m$ such that $I_{l_m} = \tilde{I}$ remains constant and $t^j_{l_m} \to \tilde{t}^j \in [0, \varepsilon/\lVert\tilde{v}\rVert]$ for $m \to + \infty$ for all $j \in \tilde{I}$. Accordingly, $x+ t^j_{l_m} v_{l_m}$ converges to $x + \tilde{t}^j \tilde{v}$ as $m \to + \infty$. 
    Since the functions $f_j$ are locally Lipschitz, there exists a common local Lipschitz constant $L\ge 0$ such that all objective functions $f_j$ are Lipschitz with constant $L$ in a neighborhood of $x+\tilde{t}^j \tilde{v}$, respectively. Due to the convergence of the sequences, we can find an index $M \ge 1$ and $\kappa \ge 0$ such that
	\begin{align} \label{eq:xi_bound}
	    {\big\|\xi^j_{l_m}\big\|}_* \le L + \kappa \quad \text{for all }m \ge M \text{ and } j \in \tilde{I}.
	\end{align}
	On the other hand, we can bound $\lVert \xi_l \rVert_* \le \lVert \xi_1 \rVert_* \le \max \lbrace \lVert \xi_1 \rVert_*, L  + \kappa \rbrace$ for all $l \ge 1$.
    For convenience, we update $L\to\max \lbrace \lVert \xi_1 \rVert_*, L + \kappa \rbrace$ for the remainder of the proof to get a uniform bound for $\lVert \xi_{l_m}^j \rVert_*$ and $\lVert \xi_l \rVert_*$ for all $m \ge M$, $j \in \tilde{I}$ and $l \ge 1$.\hfill\\
    Now, let $m \ge M$ and $j \in \tilde{I}$. Since $\xi^j_{{l_{m-1}}} \in \WW_{l_m}$ and $-\xi_{l_{m-1}} \in \conv\left(\WW_{l_{m-1}}\right) \subseteq \conv\left(\WW_{l_m}\right)$, the convex combination $(1-s)(-\xi_{l_{m-1}})+s\xi_{l_{m-1}}^j$ for $s\in (0,1)$ is in $\conv\left(\WW_{l_m}\right)$. Therefore the minimization property of $\xi_{l_m}$ yields that
	\begin{align}
        \label{eq:proof_algo_est}
		\begin{aligned}
		    &{\big\|\xi_{l_m}\big\|}_*^2 \leq{\big\|-\xi_{l_{m-1}} + s (\xi^j_{l_{m-1}} + \xi_{l_{m-1}})\big\|}_*^2\\
		    &\quad={\big\|\xi_{l_{m-1}}\big\|}_*^2-2s\,{\big\langle\xi_{l_{m-1}},\xi^j_{l_{m-1}}+ \xi_{l_{m-1}}\big\rangle}_*+s^2\,{\big\|\xi^j_{l_{m-1}} + \xi_{l_{m-1}}\big\|}_*^2\\
		    &\quad={\big\|\xi_{l_{m-1}}\big\|}_*^2 - 2 s\,{\big\langle\xi_{l_{m-1}},\xi^j_{l_{m-1}}\big\rangle}_*-2s\,{\big\|\xi_{l_{m-1}}\big\|}_*^2 + s^2\,{\big\|\xi^j_{l_{m-1}} + \xi_{l_{m-1}}\big\|}_*^2
		\end{aligned} 
	\end{align}
	for all $s \in [0,1]$. Since $j \in \tilde{I}$ we must have
	\begin{align} \label{eq:est_1}
		{\big\langle \xi_{l_{m-1}}, \xi^j_{l_{m-1}}\big\rangle}_* > -c\,{\big\|\xi_{l_{m-1}}\big\|}_*^2
	\end{align}	
	by Step 11. From inequality \eqref{eq:xi_bound} and the choice of the Lipschitz constant $L$, we can conclude that
	\begin{align} \label{eq:est_2}
		{\big\|\xi^j_{l_{m-1}} + \xi_{l_{m-1}}\big\|}_*\leq 2 L.
	\end{align}
	Combining \eqref{eq:proof_algo_est} with \eqref{eq:est_1} and \eqref{eq:est_2} yields
	\begin{align*}
		{\big\|\xi_{l_m}\big\|}_*^2 &<{\big\|\xi_{l_{m-1}}\big\|}_*^2 + 2 s c\,{\big\|\xi_{l_{m-1}}\big\|}_*^2 - 2 s\,{\big\|\xi_{l_{m-1}}\big\|}_*^2+4s^2 L^2\\
		&={\big\|\xi_{l_{m-1}}\big\|}_*^2 - 2 s (1-c)\,{\big\|\xi_{l_{m-1}} \big\|}_*^2 + 4 s^2 L^2.
	\end{align*}   
	Now, we choose $s \coloneqq (1-c)\|\xi_{l_{m-1}}\|_*^2/(4L^2)$. Since $1 - c \in (0,1)$ and $\lVert v_{l_{m-1}} \rVert_* \leq L$ we have $s \in(0,1)$. Thus, we obtain
	\begin{align*}
		{\big\|\xi_{l_m}\big\|}_*^2 &<{\big\|\xi_{l_{m-1}}\big\|}_*^2-\frac{2(1-c)^2}{4 L^2}\,{\big\|\xi_{l_{m-1}}\big\|}_*^4 + \frac{(1-c)^2}{4 L^2}\,{\big\|\xi_{l_{m-1}}\big\|}_*^4 \\
		&= \left( 1 - \frac{(1-c)^2}{4 L^2}\,{\big\|\xi_{l_{m-1}}\big\|}_*^2 \right){\big\|\xi_{l_{m-1}}\big\|}_*^2.
	\end{align*}    
	We have assumed that Algorithm \ref{algo:descent_direction} does not terminate. Therefore, we must have $\lVert \xi_{l_{m-1}} \rVert_* > \delta$, which implies
	\begin{align*}
		{\big\|\xi_{l_m}\big\|}_*^2 < \underbrace{\left( 1 - \left( \frac{1-c}{2 L} \,\delta \right)^2 \right)}_{\eqqcolon r} {\big\|\xi_{l_{m-1}}\big\|}_*^2.
	\end{align*}    	
	Note that we have $\delta < \lVert \xi_{l_m} \rVert_* \leq L$ for all $l \in \N$, so $r \in (0,1)$. Additionally, $r$ does not depend on $l_m$, so we have
	\begin{align*}
		{\big\|\xi_{l_m}\big\|}_*^2< r\,{\big\|\xi_{l_{m-1}}\big\|}_*^2< r^2\,{\big\|\xi_{l_{m-2}}\big\|}_*^2<\ldots< r^{m-1}\,{\big\|\xi_{l_1}\big\|}_*^2 \leq r^m L^2.
	\end{align*}	
	In particular, there is some $m$ such that $\lVert \xi_{l_m} \rVert_* \leq \delta$, which is a contradiction. 
\end{proof}
	
\begin{remark}
	The proof of Theorem~{\em \ref{thm:descent_dir_conv}} shows that for convergence of Algorithm~{\em\ref{algo:descent_direction}}, it would be sufficient to consider only a single $j \in I_l$ in Step 11. Similarly, for the initial approximation $\WW_1$, a single element of $\partial_\varepsilon f_i(x)$ for any $i \in \{1,\dots,k\}$ would be enough. A modification of either step can potentially reduce the number of executions of Step 11 (i.e., Algorithm {\em\ref{algo:new_subgradient}}) in Algorithm {\em\ref{algo:descent_direction}} in case the $\varepsilon$-subdifferentials of multiple objective functions are similar. However, we will forgo these modifications and leave Algorithm {\em\ref{algo:descent_direction}} as it is, since both modifications also introduce a bias towards certain objective functions, which we want to avoid.\hfill$\Diamond$
\end{remark}

\subsection{A descent method for nonsmooth MOPs}
\label{subsec:descent_method_nonsmooth_MOPs}

Building on Algorithm \ref{algo:descent_direction}, it is now straightforward to construct the descent method for locally Lipschitz continuous MOPs given in Algorithm~\ref{algo:nonsmooth_descent_method}.
\begin{algorithm} 
	\caption{(Nonsmooth descent method)}
	\label{algo:nonsmooth_descent_method}
	\begin{algorithmic}[1] 
		\REQUIRE Initial point $x_1 \in \H$, parameters for stopping criterion $\overline{\delta}, \overline{\varepsilon} \ge 0$, tolerance sequences $(\delta_j)_j, (\varepsilon_j)_j \subseteq \R_{>0}$, Armijo parameters $c \in (0,1)$, $t_0 > 0$.
		\FOR{$j = 1,2,\dots$}
			\STATE Compute a descent direction $v_j$ via Algorithm \ref{algo:descent_direction} with inputs $(x_j, \varepsilon_j, \delta_j, c)$.
			\STATE Use backtracking line search to determine
				\begin{align*}
					\bar{s}
					=\inf\big\{ s \in \N \cup\big\{ 0 \}\,\big|&~f_i(x_j + 2^{-s} t_0 v_j) \leq f_i(x_j) - 2^{-s} ct_0\,\|v_j\|^2\\
				&~\text{for all }i \in \{1,\dots,k\} \big\}
				\end{align*}

				and set $\bar{t} = \max\{ 2^{-\bar{s}} t_0,\varepsilon_j/\lVert v_j \rVert \}$.
            \IF{$\lVert v_j \rVert \le \overline{\delta}$ and $\varepsilon_j \le \overline{\varepsilon}$}
                \RETURN $(\overline{\varepsilon}, \overline{\delta})$-critical point $x_j$
            \ELSE
			    \STATE Set $x_{j+1} = x_j + \bar{t} v_j$.
            \ENDIF
		\ENDFOR
	\end{algorithmic} 
\end{algorithm}	
In Step 3, the classical Armijo backtracking line search is used (cf.~\cite{FS2000}) for the sake of simplicity. Note that it is well-defined due to Step 7 in Algorithm~\ref{algo:descent_direction}.
	
Clearly, the stopping condition matches the Definition \ref{def:eps_delat_crit} of the current iterate being $(\overline{\varepsilon},\overline{\delta})$-critical exactly. Thus, when Algorithm~\ref{algo:nonsmooth_descent_method} terminates, it will in fact return an $(\overline{\varepsilon},\overline{\delta})$-critical point.
We state a convergence as well as a termination result for Algorithm \ref{algo:nonsmooth_descent_method}. 
First off, in Theorem \ref{thm:convergence_proof_nonsmooth_descent_method}, we address the case, where the tolerances $\overline{\varepsilon}$ and $ \overline{\delta}$ are both set to $0$. The theorem states that we converge (in the sense of subsequences) to Pareto critical points in the limit. Then, in Theorem \ref{thm:convergence_proof_nonsmooth_descent_method_2} we show that the algorithm is capable of finding $(\overline{\varepsilon},\overline{\delta})$-critical points, for generalized parameter settings. 

\begin{theorem}
    \label{thm:convergence_proof_nonsmooth_descent_method}
    For $i = 1,\dots, k$ let $f_i\colon\H\to\R$ be locally Lipschitz. We suppose that $x_1$ is an element in $\H$ and $(\delta_j)_j, (\varepsilon_j)_j \subseteq \R_{>0}$ be two sequences with
    \begin{align*}
        \delta_j \to 0,\quad\varepsilon_j \to 0\quad\text{ and }\quad\sum_{j=1}^{\infty} \varepsilon_j\delta_j=\infty. 
    \end{align*}
    Let further $\overline{\varepsilon} = \overline{\delta} = 0$, $c\in (0,1)$ and $t_0>0$.
    Assume Algorithm {\em\ref{algo:nonsmooth_descent_method}} does not converge after finitely many steps.
    Let $(x_j)_j$ be the sequence generated by Algorithm~{\em\ref{algo:nonsmooth_descent_method}} with inputs $(x_1, \overline{\delta}, \overline{\varepsilon}, (\delta_j)_j, (\varepsilon_j)_j,c,t_0)$. Then, we have:
	\begin{itemize}
	    \item[\rm a)] Every accumulation point of $(x_j)_j$ is Pareto critical.
		\item[\rm b)] If one $f_i$ is bounded from below then $(x_j)_j$ possesses a subsequence $(x_{j_l})_{l}$ such that $\lVert v_{j_l}\rVert \to 0$ as $l \to \infty$.
    \end{itemize}
\end{theorem}
	
\begin{proof}
    \begin{itemize}
        \item [a)] In the following proof we choose appropriate subsequences of $(x_j)_j$ multiple times. We will do this without relabeling the sequence and only comment when doing so. Let $x^*$ be an accumulation point of $(x_j)_j$. Then, there exists a subsequence (no relabeling) with $x_j \to x^*$ as $j \to \infty$.\hfill\\ 
	   First we show that $\lVert v_j\rVert \le \delta_j$ is true for infinitely many $j \ge 1$. In each iteration of Algorithm \ref{algo:nonsmooth_descent_method}, we use Algorithm \ref{algo:descent_direction}. Therefore at least one of the stopping criteria of Algorithm \ref{algo:descent_direction} is met infinitely many times. Assume the stopping criteria $\lVert v_l \rVert < \delta$ in Step 4 of Algorithm \ref{algo:descent_direction} (where $\lVert v_l\rVert = \lVert \xi_l\rVert_*$) is only met finitely many times. Then, there exists $l \ge 1$ such that for all $j\ge l$ it holds that
	   \begin{align}
	       \label{eq:stopping_criterions}
	    f_i (x_{j+1} ) \le f_i(x_j) - c \varepsilon_j \lVert v_j \rVert\text{ for all }i = 1, \ldots, k \quad\text{and}\quad\lVert v_j \rVert &> \delta_j.
	   \end{align}
	   The first inequality follows from the active stopping criterion in Step 7 of Algorithm \ref{algo:descent_direction} and the way the backtracking rule in Step 3 of Algorithm \ref{algo:nonsmooth_descent_method} is defined. We show that these inequalities lead to a contradiction. Let $i \in \{1, \dots, k\}$ and $J \ge l$. Then, we have
	   \begin{align}\label{eq:descent_contradiction}
	    \begin{split}
	     f_i(x_{J+1}) - f_i(x_{l}) & = \sum_{j = l}^{J} f_i(x_{j+1}) - f_i(x_j) \le \sum_{j = l}^{J} -c \varepsilon_j \lVert v_j \rVert\\
        &< -c \sum_{j = l}^{J} \varepsilon_j\delta_j.
	     \end{split}
	   \end{align}

	   We know by the assumptions on $(\delta_j)_j$ and $(\varepsilon_j)_j$ that the last series diverges to $\infty$. 
        Accordingly, the sequential continuity of $f_i$ yields that
            \begin{align*}
    	       f_i(x^*) - f_i(x_{l}) = -c \lim_{J \to +\infty} \sum_{j = l}^{J-1} \varepsilon_j\delta_j = -\infty.
           \end{align*}
        which is a contradiction as the difference on the left hand side is finite.

	  Therefore, $\lVert v_j\rVert \le \delta_j$ holds for infinitely many $j \ge 1$. 
   This means, we can choose an appropriate subsequence of $(x_j)_j$ (no relabeling) such that
	   \begin{align*}
	       x_j\to x^*\text{ as }j \to \infty \quad \text{ and}\quad\lVert v_j\rVert<\delta_j\text{  for all } j \ge J.
       \end{align*}       
	   By Theorem \ref{thm:subdiff_eps_subdiff_multi} it follows that $0 \in \conv(\cup_{i=1}^k \partial f_i(x^*))$. Hence $x^*$ is Pareto critical.
	   \item [\rm b)] The proof follows from inequalities \eqref{eq:stopping_criterions} and \eqref{eq:descent_contradiction} and the fact that $\lVert v_j \rVert \le \delta_j$ has to hold for infinitely many $j \ge 1$ if $f_i$ is bounded from below.
    \end{itemize}    
\end{proof}

In practice, we will of course rely on Algorithm \ref{algo:nonsmooth_descent_method} terminating after a finite number of iterations due to the stopping criterion for tolerances $\overline{\varepsilon}, \overline{\delta} > 0$ instead of generating infinite sequences of iterates.
The following theorem states that the algorithm will in fact terminate after a finite number of iterations, e.g., if the sequences $(\varepsilon_j)_j$ and $(\delta_j)_j$ are chosen as certain constants.

\begin{theorem}
    \label{thm:convergence_proof_nonsmooth_descent_method_2}
    For $i = 1,\dots, k$ let $f_i\colon\H\to\R$ be locally Lipschitz continuous. We suppose that $x_1$ is an element in $\H$ and set $\overline{\varepsilon}, \overline{\delta} > 0$. Let $(\delta_j)_j, (\varepsilon_j)_j \subseteq \R_{>0}$ be constant sequences with $\delta_j = \overline{\delta}$, $\varepsilon_j = \overline{\varepsilon}$, $c\in (0,1)$ and $t_0>0$. Let $(x_j)_j$ be the sequence generated by Algorithm~{\em\ref{algo:nonsmooth_descent_method}} with inputs $(x_1, \overline{\delta}, \overline{\varepsilon}, (\delta_j)_j, (\varepsilon_j)_j,c,t_0)$ and assume that one objective function $f_i$ is bounded from below. Then Algorithm~{\em\ref{algo:nonsmooth_descent_method}} returns an $(\overline{\varepsilon}, \overline{\delta})$-critical point after finitely many iterations.
\end{theorem}

\begin{proof}
    Assume Algorithm \ref{algo:nonsmooth_descent_method} does not terminate after finitely many steps and produces an infinite sequence $(x_j)_j$. Since the condition $\varepsilon_j \le \overline{\varepsilon}$ is fulfilled in every iteration of Algorithm \ref{algo:nonsmooth_descent_method}, we show that $\lVert v_j \rVert \le \overline{\delta}$ has to hold for one $j \ge 1$. Then Algorithm \ref{algo:nonsmooth_descent_method} stops since the condition $\varepsilon_j \le \overline{\varepsilon}$ is fulfilled in every step. Again one of the stopping criterion of Algorithm \ref{algo:descent_direction} has to be fulfilled infinitely many times. If $\lVert v_j \rVert \le \overline{\delta}$ in Step 4 of Algorithm \ref{algo:descent_direction} is fulfilled then also Algorithm \ref{algo:nonsmooth_descent_method} stops. If this is not the case then Algorithm \ref{algo:descent_direction} only stops due to the stopping condition in Step 8 and we conclude that for all $j \ge 1$ it holds that
    \begin{align}
	  \label{eq:stopping_criterions_2}
           f_i (x_{j+1}) \le f_i(x_j) - c\,\overline{\varepsilon}\lVert v_j \rVert\text{ for all }i = 1, \ldots, k \quad\text{and}\quad\lVert v_j \rVert &> \overline{\delta}.
	  \end{align}
    Combining these inequalities we have for all $j \ge 1$
    \begin{align}
        \label{eq:stopping_criterions_3}
           f_i (x_{j+1}) \le f_i(x_j) - c\, \overline{\varepsilon}\overline{\delta} \text{ for all }i = 1, \ldots, k.
    \end{align}
    This leads to a contradiction. Fix $i \in \{1,\dots ,k\}$ such that $f_i$ is bounded from below. Then for all $J \ge 1$ we have by \eqref{eq:stopping_criterions_3}
    \begin{align}
        \label{eq:telescope_sum_2}
        f_i(x_{J+1}) - f_i(x_1) = \sum_{j = 1}^{J} f_i(x_{j+1}) - f_i(x_j) \le \sum_{j = 1}^{J} -c\, \overline{\varepsilon} \overline{\delta} = - J c\, \overline{\varepsilon}\overline{\delta}.
    \end{align}
    Since the right-hand side of \eqref{eq:telescope_sum_2} diverges to $-\infty$ for $J \to \infty$, we arrive at a contradiction given that $f_i$ is bounded from below. 
    \end{proof}

\begin{remark}
    The choice of the tolerance sequences $(\delta_j)_j$ and $(\varepsilon_j)_{j}$ in Theorem~\ref{thm:convergence_proof_nonsmooth_descent_method_2} can be further relaxed. We are not forced to use constant sequences $\delta_j = \overline{\delta}$ and $\varepsilon_j = \overline{\varepsilon}$. Instead, we could choose arbitrary sequences with $\delta_j \in (0, \overline{\delta}]$ and $\varepsilon_j \in (0, \overline{\varepsilon}]$ that satisfy the condition $\sum_{j=1}^{\infty} \delta_j \varepsilon_j = \infty$ similar to the requirements of Theorem~\ref{thm:convergence_proof_nonsmooth_descent_method}. This could be further relaxed to arbitrary positive sequences $\delta_j > 0$ and $\varepsilon_j > 0$ provided that they remain bounded by $\overline{\delta}$ and $\overline{\varepsilon}$ for almost all iterations and that they also satisfy the summability property  $\sum_{j=1}^{\infty} \delta_j \varepsilon_j =\infty$. The proof in these settings follows analogously to the proof of Theorem~\ref{thm:convergence_proof_nonsmooth_descent_method_2}.\hfill$\Diamond$
\end{remark}

\section{Application in bicriterial optimal control of an obstacle problem}
\label{sec:Application}

In this section, we examine the behavior of Algorithm \ref{algo:nonsmooth_descent_method} applied to a classic, nonsmooth obstacle-constrained optimal control problem -- see, e.g., \cite[Section 6]{KinderlehrerStampacchia:1980:1} -- on the two-dimensional domain $\Omega \coloneqq (-1,1)^2$ for two objective functions.

The forward problem, i.e., the constraint in the optimal control problem, can be interpreted as the problem of finding a displacement $y\colon \Omega \to \R$ of a clamped membrane under external, distributed vertical forces $u\colon \Omega \to \R$ (assuming small displacements with linear response) with a rigid obstacle, described by $\psi\colon\Omega\to\R$, limiting the vertical displacement to $y\leq \psi$.

This constrained problem can be equivalently formulated as a convex energy minimization problem or via the corresponding partial differential variational inequality, and it is well understood.
Most importantly, the control-to-state operator is known to be well-defined, Lipschitz continuous and Hadamard- but generally not Fréchet-differentiable everywhere, cf.
\cite{LS67, Har77, Mig76}. There is also extensive literature on computational aspects for obstacle constrained dynamics, including efficient solvers, cf. \cite{WW10,MRW15,GK09,GSS09}.

Various aspects of optimal control problems with the obstacle constraint have previously been considered in a broad range of publications (e.g., \cite{Mig76, Wac14 ,Christof20}), but, to the best of our knowledge, obstacle-constrained optimization problems have not been considered in the context of infinite-dimensional multiobjective optimization (though their discretizations have been dealt with in finite-dimensional, nonsmooth multiobjective optimization, c.f.~\cite{MM1993}). Due to the nonlinearity of the control-to-state operator, these problems are generally nonconvex and nonsmooth. However, (varying notions of) subdifferentials of the control-to-state operator have been characterized in \cite{RW19}, and \cite[Theorem 5.7]{RU19} shows how to compute an element of the Clarke subdifferential of control reduced optimal control of the obstacle problem -- which is what we require in order to employ our common-descent method. Note that this exact technique for computing subderivatives was applied in scalar optimal control of obstacle-constrained problems using an inexact bundle method in function space, see \cite{hertlein2019inexact, hertlein2022inexact}.

\subsection{Problem description}
\label{subsec:prob_desc}

The domain we consider is the two dimensional square $\Omega = (-1,1)^2 \subseteq \R^2$ with an obstacle described by a function $\psi \in H^1(\Omega)$ (to be specified later), yielding the set of admissible displacements
\begin{align*}
    K \coloneqq \left\lbrace y \in V\coloneqq \Vspace : y \le \psi \text{ a.e. on } \Omega \right\rbrace,
\end{align*}
which is guaranteed to be nonempty by choosing $\psi$ appropriately. The variational inequality formulation of the constraining obstacle problem for a fixed, distributed external load $f\in V'\coloneqq \Vdual$ amounts to the finding $y \in K$, such that
\begin{align}
    \label{obstacleVI}
    {\langle Ay - f, v - y \rangle}_{V',V} \ge 0 \, \text{ for all } v \in K
    .
\end{align}
Here, $A\colon V \to V'$ is a linear, continuous and coercive partial differential operator (we will be using the weak form of $A = -\Delta$ in the following), and $\langle \cdot\,, \cdot \rangle_{V',V}$ denotes the dual pairing. In the optimal control problem, we consider the control space $U\coloneqq L^2(\Omega)$ with the standard $U \hookrightarrow V'$ Gelfand-type embedding to let $u \in U$ assume the role of $f$ in \eqref{obstacleVI}.

Given a desired state and reference control $y_d \in H\coloneqq L^2(\Omega)$, $u_d \in U$, we then fix the two cost functionals to obtain the optimal control problem
\begin{equation}
    \label{prob:obstacle}
    \begin{aligned}
        \min_{(y,u)\in K\times U}&~\frac{1}{2}
        \begin{pmatrix}
			\lVert y-y_d\rVert_H^2\\[1mm] 
            C\,\lVert u - u_d \rVert_U^2 
		\end{pmatrix},&\\
        \text{ s.t. } &{\langle Ay - u, v - y \rangle}_{V',V} \ge 0 \, \text{ for all } v \in K,
    \end{aligned}
\end{equation}
with a hyperparameter $C>0$. Note that $C$ is essentially introduced in order to scale the axes in the plots of the Pareto fronts, so that they are easier to interpret. Introducing and tuning the parameter $C$ can be interpreted as preconditioning of the problem.

Problem \eqref{prob:obstacle} is an optimal control problem and clearly a constrained problem. To make it fit into the realm of unconstrained optimization, which we have formulated the algorithms in this paper for, we simply make use of the existence of the Hadamard-differentiable solution operator of the obstacle problem $S\colon U\to K\subseteq V$ mapping a control $u$ to the solution $y=S(u)$ of the constraining variational inequality of \eqref{prob:obstacle} to obtain the equivalent control-reduced multiobjective optimization problem 
\begin{equation}
\begin{aligned}
    \min_{u \in U}&~\frac{1}{2} 
        \begin{pmatrix}
			\lVert S(u) - y_d \rVert_H^2  \\[1mm] 
        C\,\lVert u - u_d \rVert_U^2 
		\end{pmatrix}.&
    \label{prob:obstacle_red}
\end{aligned}
\end{equation}

Using the direct method of variational calculus, one can easily show, that the weighted-sum-scalarized problems corresponding to this problem possess solutions, and hence the Pareto set and the Pareto front of this problem are nonempty.
What remains to be fixed in the remainder is the choice of the algorithmic parameters, the desired states and controls $y_d,u_d$ and the specific obstacle $\psi$. 

We describe the choice of the free parameters in the following subsection. In all cases, we ensure that our problem configuration in fact captures the nonsmooth behavior of the problem.
As mentioned above, the nonsmoothness of  the problem is introduced by the solution operator.
More specifically, the points of non-Fréchet-differentiability are precisely those of so called \emph{weak contact}, i.e., where the control corresponds to a state that is in contact with the obstacle, but where there are no normal forces actively preventing penetration on a sufficiently large area (in the sense of Sobolev capacities). Such configurations of "coincidental" contact are exactly those, where the problem transfers from a free Poisson problem to a full constrained problem.

\subsection{Computational procedure and joint parameters}
\label{subsec:comp_proc}

The goal of our numerical procedure is to find an approximate representation of the Pareto front and Pareto set of the obstacle-constrained optimal control problem \eqref{prob:obstacle_red}. To this end, we apply Algorithm \ref{algo:nonsmooth_descent_method} starting from a number of varying initial values. As shown in Theorem \ref{thm:convergence_proof_nonsmooth_descent_method_2}, for each initial value, Algorithm \ref{algo:nonsmooth_descent_method} terminates at an $(\overline{\varepsilon}, \overline{\delta})$-critical point after finitely many steps. As the terminal iterate of the algorithm typically varies with varying initial guesses, we obtain a representation of the Pareto front and the Pareto set of \eqref{prob:obstacle_red} by $(\overline{\varepsilon}, \overline{\delta})$-critical points. We chose the different initial controls $u_0 \in U$ constant on the entire domain. Specifically, we apply the algorithm for constant initial controls for all values $u_0 \equiv \hat u\in \{ 1, 2, \dots, 8\}$ and for all mesh discretizations $h_{\max} \in \{0.2, 0.1, 0.05, 0.02\}$.
    
For all experiments, we fix the scaling parameter $C=1.5e$--$2$ and the hyper\-parameters $\overline{\varepsilon} = \overline{\delta} = 1e$--$4$, $c = 1e$--$1$ and the constant sequences $(\varepsilon_i)_i \equiv \overline{\varepsilon}, (\delta_i)_i \equiv \overline{\delta}$. Further, we set $y_d \equiv 2$ and $u_d \equiv 0$. This choice yields a setting where the first cost functional improves when the state is pushed upwards towards the desired state, while the second objective is optimal for vanishing controls, leading to a setting where optimal compromises can be expected to achieve some upwards deformation of the state using controls "efficiently". This suggests that contact should be established in optimal compromises, but no additional forces are to be applied, leading to a nonsmooth weak-contact situation in the optimal compromise. 
We therefore expect the algorithm to have to deal with increasing nonsmoothness over the course of the run.

Note that at this point, only the obstacle remains to be fixed in each of the examples. We will specify the obstacles we use in the experiment runs in Subsections~\ref{subsub: Configuration 1 : constant obstacle} and \ref{subsub: Configuration 2 : piecewise constant obstacle}, respectively.
 
\subsection{Implementation details}
\label{subsec:imp_det}

We discretize the optimal control problem using Lagrangian $P1$ finite elements on a triangulation of $\Omega$ supplied by MATLAB's PDE-toolbox with a predetermined target maximum element edge length $h_{\text{max}}$ (which is typically only violated by fractions of a percent) and nodally interpolate the obstacle $\psi$ to essentially enforce the nonpenetration constraint nodally. The discretizations of $\Omega$ we use are those corresponding to $h_{\max} \in \{0.2, 0.1, 0.05, 0.02\}$. Additionally, we compute a reference solution $u^*_{\textrm{ref}}$ for $h_{\max}=0.01$ to emulate the exact solution in order to investigate convergence of the solutions for finer meshes. The number of finite elements corresponding to each mesh discretization can be seen in Table~\ref{tab:FEM_sizes}, ranging from $135$ to $45\,857$ elements. 
\begin{table}
    \centering
    \begin{tabular}{!{\vrule width 1.2pt}c|c|c|c|c|c!{\vrule width 1.2pt}}
        \noalign{\hrule height 1.2pt}
        $h_{\max}$ & $0.2$ & $0.1$ & $0.05$ & $0.02$ & $0.01$ \\ \hline
        \# FEM & 135 & 494 & 1\,909 & 11\,682 & 45\,857 \\ \noalign{\hrule height 1.2pt}
    \end{tabular}
    \caption{Number of finite elements for different maximum edge lengths $h_{\max}$.}
    \label{tab:FEM_sizes}
\end{table}
 
The control-to-state operator is implemented using an active-set strategy applied to the equivalent energy minimization formulation of the obstacle problem and the subderivatives are obtained based on the discretized analogue of the adjoint-based computations in \cite[Theorem 5.7]{RU19}, where the discrete approximation to the adjoint state is computed using MATLAB's \emph{mldivide} routine to solve the corresponding linear system. Our implementations of Algorithms \ref{algo:new_subgradient}-\ref{algo:nonsmooth_descent_method} is also in MATLAB. The preconditioner that maps generalized subderivatives to primal objects, i.e., Riesz's operator (in, e.g., Lemma~\ref{lem:step_size_bound}), is chosen as the canonical $\Hspace$-Riesz operator.

\subsection{Numerical results}

In this subsection, we present the numerical results obtained by Algorithm~\ref{algo:nonsmooth_descent_method} for the optimal control problem described in Subsection~\ref{subsec:prob_desc}. The settings of the parameters for Algorithm~\ref{algo:nonsmooth_descent_method} are specified in Subsection~\ref{subsec:comp_proc}, while the implementation details to handle the PDE-constraints are described in Subsection~\ref{subsec:imp_det}. To conduct the experiments, we only have to choose the shape of the obstacle $\psi$, which we do in two example instances below. We consider a constant obstacle in Subsection~\ref{subsub: Configuration 1 : constant obstacle} and a more involved example in Subsection~\ref{subsub: Configuration 2 : piecewise constant obstacle}. Further, in Subsection~\ref{subsub:experiment_size_exp_subdiff}, we analyze the size of the approximated Goldstein $\varepsilon$-subdifferential, which is computed in every iteration of Algorithm~\ref{algo:nonsmooth_descent_method} using Algorithm~\ref{algo:descent_direction}, in order to investigate the behaviour of our algorithm.

\begin{figure}[h]  
    \begin{center}
    \begin{subfigure}[t]{.28\textwidth}
        \centering
        \includegraphics[width = \linewidth]{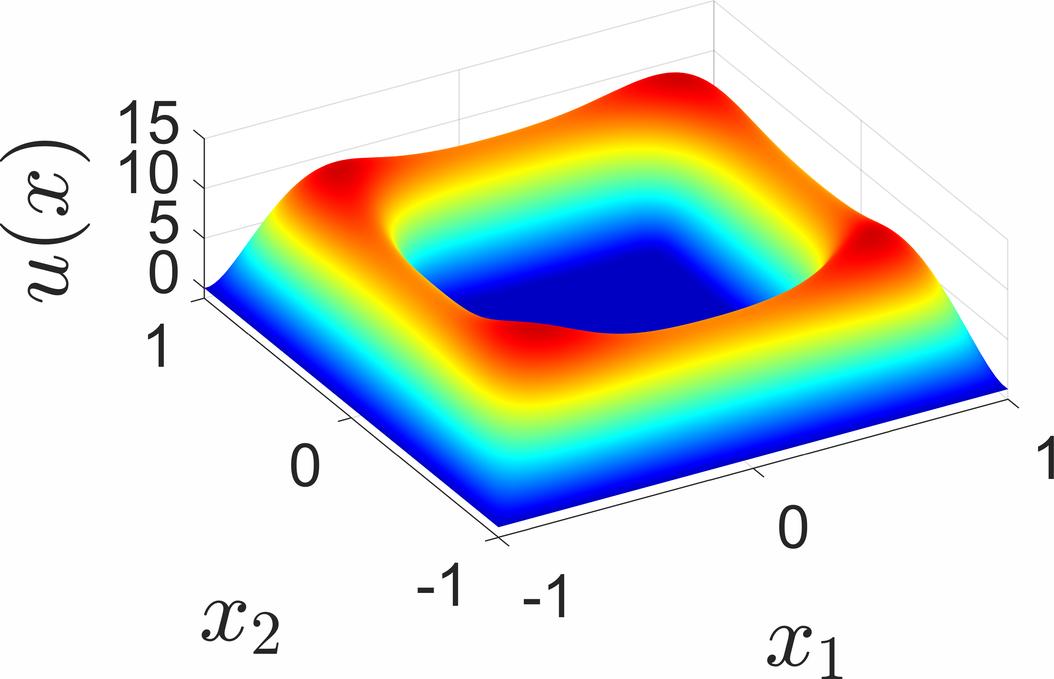}
        \caption{Control}
        \label{subfig:const_control}
    \end{subfigure}
    \hfill
    \begin{subfigure}[t]{.28\textwidth}
        \centering
        \includegraphics[width=\linewidth]{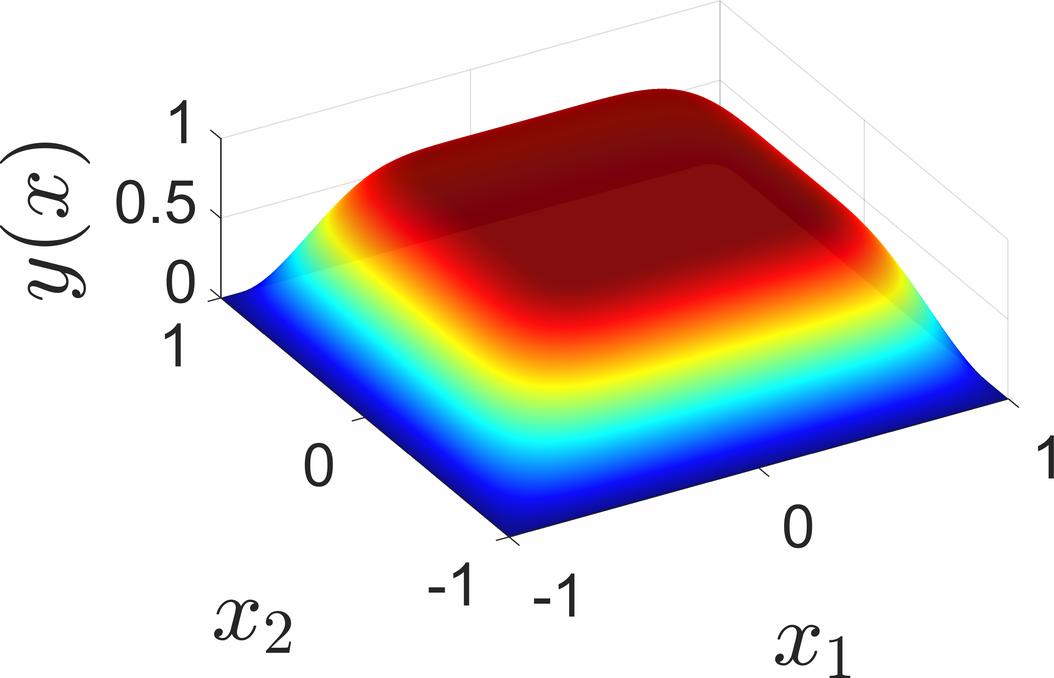}
        \caption{State}
        \label{subfig:const_state}
    \end{subfigure}
    \hfill
    \begin{subfigure}[t]{.28\textwidth}
        \centering
        \includegraphics[width = \linewidth]{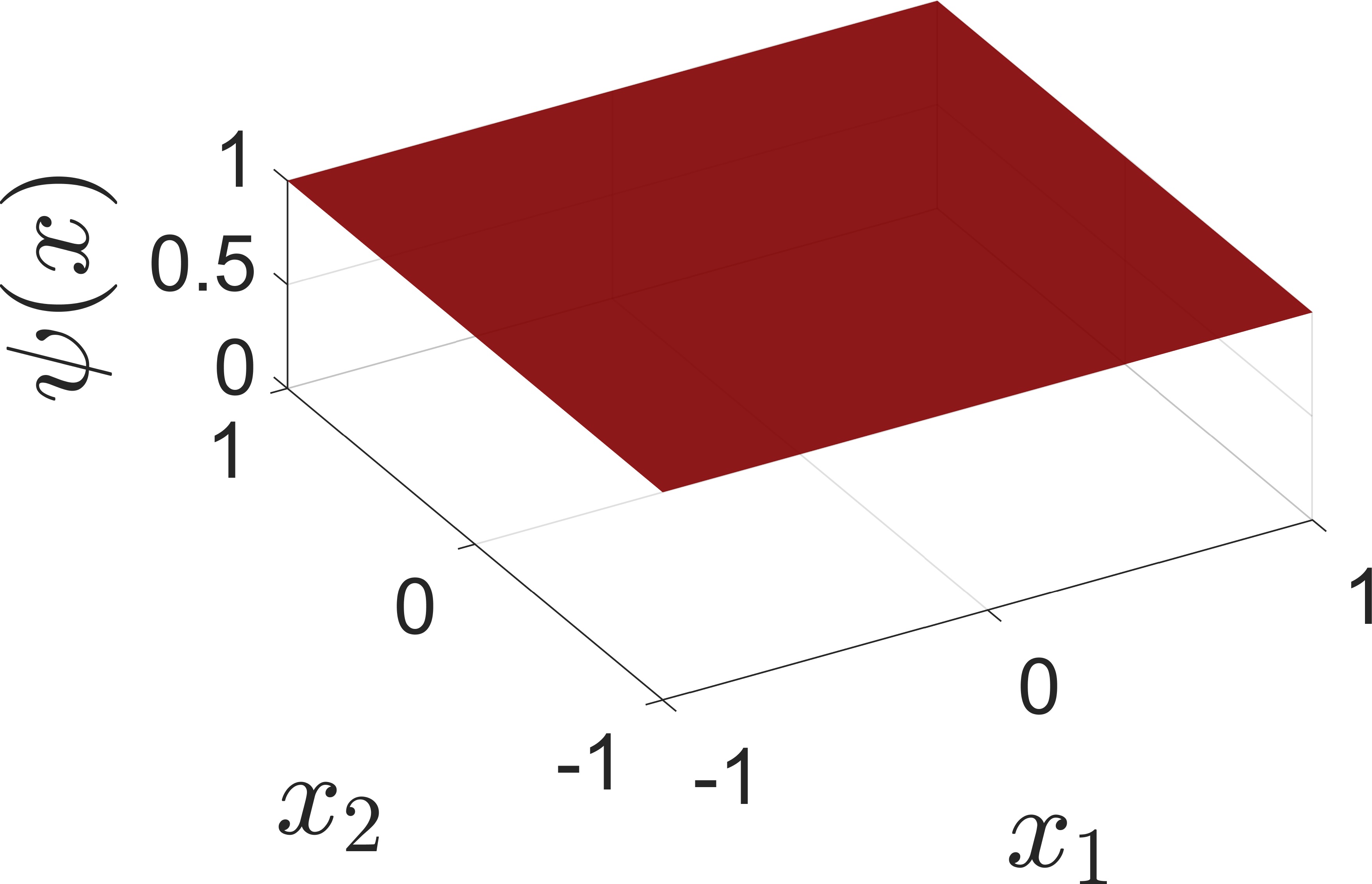}
        \caption{Obstacle}
        \label{subfig:const_obstacle}
    \end{subfigure}
    \caption{A Pareto optimal control computed with Algorithm~\ref{algo:nonsmooth_descent_method} for mesh size $h_{\max} = 0.02$, initial control $u_0 \equiv 8$ and the constant obstacle $\psi \equiv 1$.}
    \label{fig:const_pareto_sol}    
    \end{center}
\end{figure}
    
\subsubsection{Configuration 1: Constant obstacle}
\label{subsub: Configuration 1 : constant obstacle}

For the first example configuration, we set $\psi \equiv 1$. Since the desired state is $y_d \equiv 2$, the minimization of $J_1(u) = \nicefrac{1}{2}\,\lVert S(u) - y_d \rVert_H^2$ is expected to lead to configurations with contact $y(x) = \psi(x)$ for some points $x \in \Omega$. On the other hand, the second objective function $J_2(u) = \nicefrac{C}{2}\,\lVert u - u_d \rVert_U^2$, with $u_d \equiv 0$, enforces vanishing controls. We end up in a scenario with conflicting objective functions, with solutions drawn to the obstacle by one objective. An (approximate) optimal compromise in this conflicting setting can be seen in Figure~\ref{fig:const_pareto_sol}. 
Subfigure~\ref{subfig:const_control} shows the optimal control $u$ computed over $2\,697$ iterations. The corresponding state $y$ is shown in Subfigure~\ref{subfig:const_state} with the obstacle $\psi$ in Subfigure~\ref{subfig:const_obstacle}. All solutions obtained by Algorithm~\ref{algo:nonsmooth_descent_method} for the different meshes and initial states share similar features. In the middle of the domain, there is an area of contact, i.e., a region with $y(x) = \psi(x)$. In this area the control $u(x)$ vanishes. This is intuitive, since increasing the control at a point with contact only increases the objective function value of $\nicefrac{C}{2}\,\lVert u - u_d \rVert^2_U$ without decreasing the objective function value of $\nicefrac{1}{2}\,\lVert S(u) - y_d \Vert^2_H$. The size of the area of contact is influenced by the magnitude of the initial control $u_0$. For larger control values, we observe a larger area of contact in the solution, while for smaller values, the size of the area of contact is smaller. If we start with small initial control (e.g., $u_0 \equiv 1$), we get solutions with no contact at all, i.e., solutions where the obstacle problem reduces to Poisson's equation and the obstacle $\psi$ can be ignored.
    
\begin{figure}  
\begin{center}
    \begin{subfigure}[t]{.401\textwidth}
        \centering
        \includegraphics[width = \linewidth]{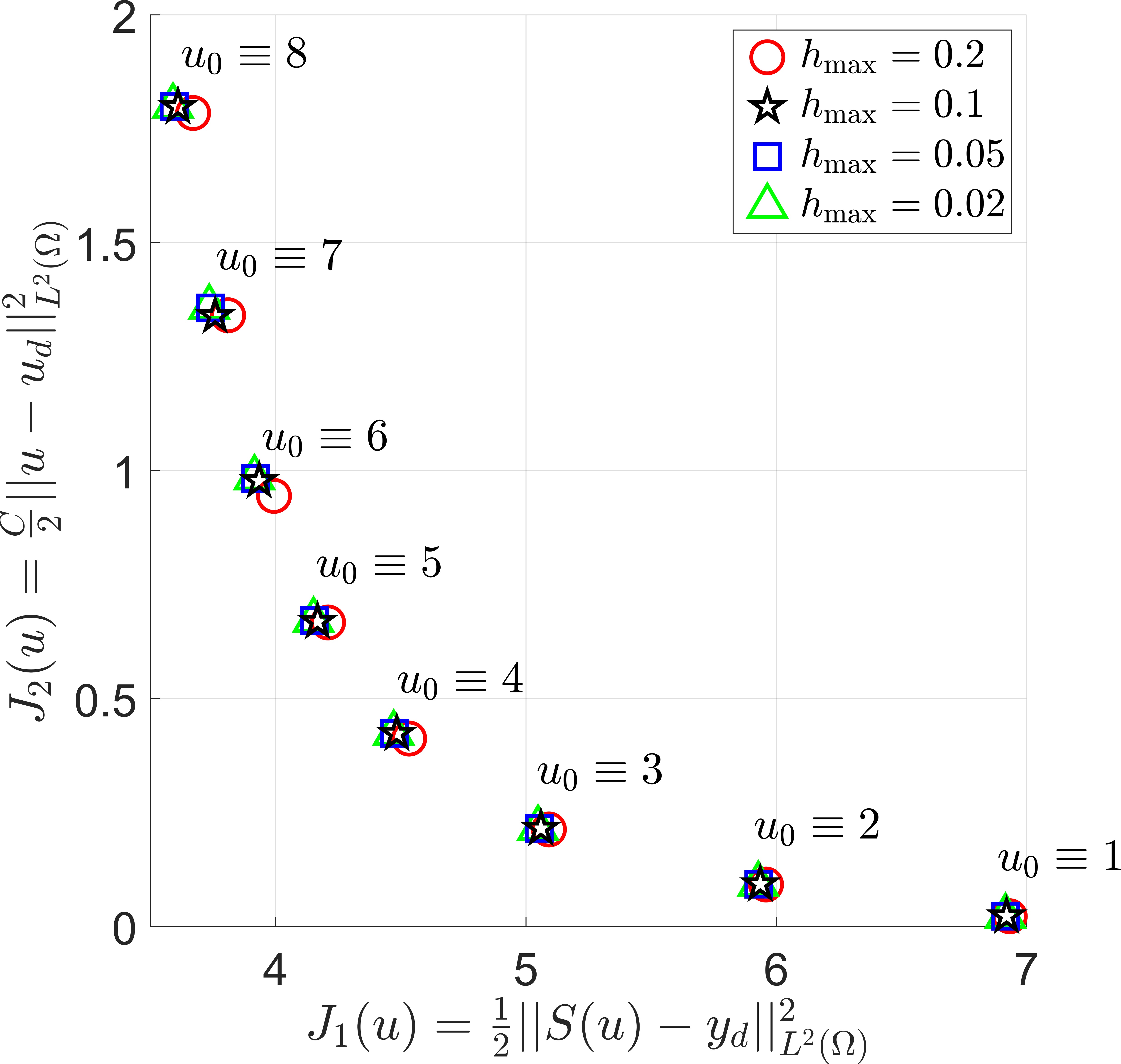}
        \caption{Pareto fronts for different mesh discretizations.}
        \label{subfig:const_parefronts}
    \end{subfigure}
    \hspace{10mm}
    \begin{subfigure}[t]{.401\textwidth}
        \centering
        \includegraphics[width=\linewidth]{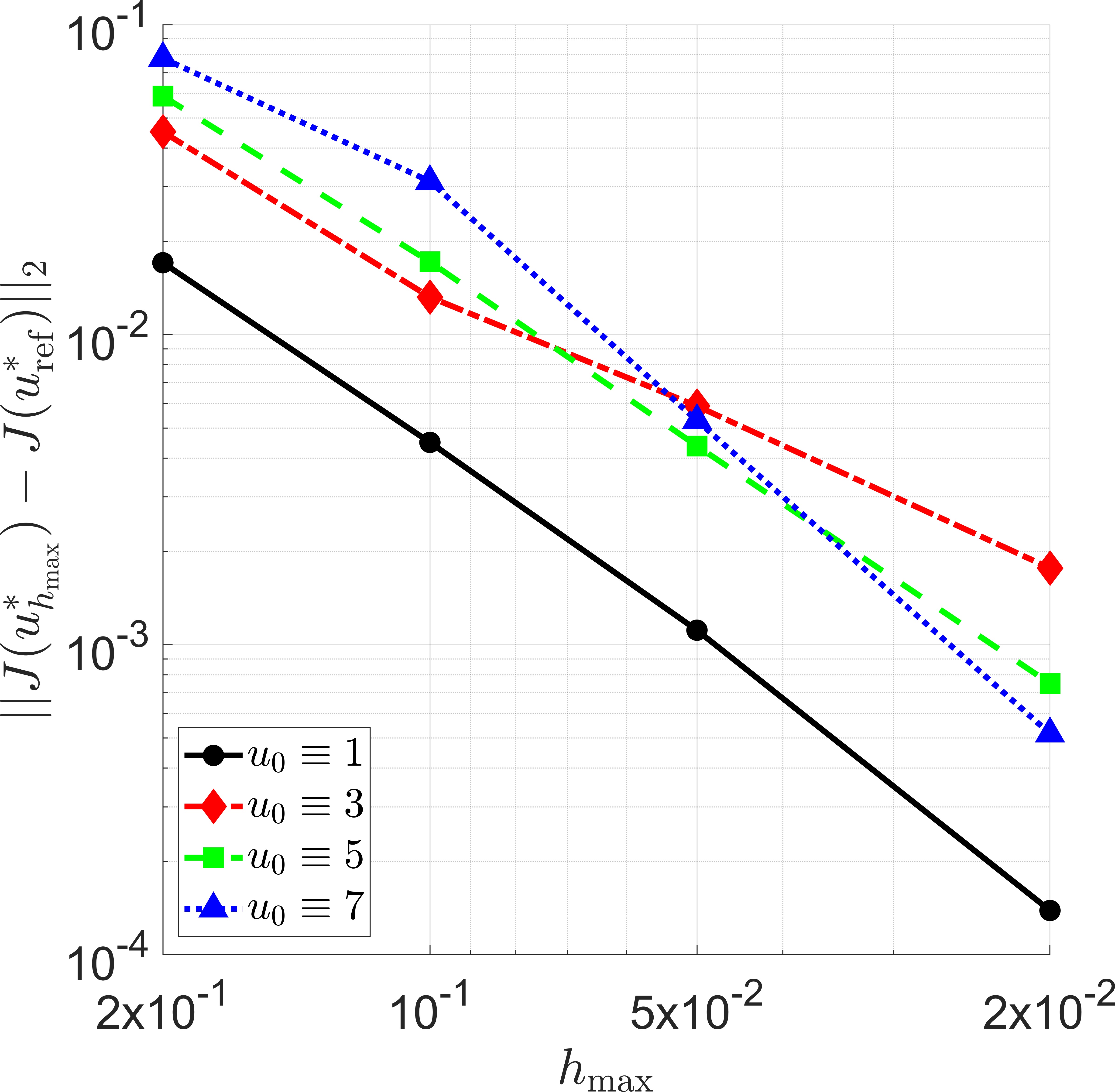}
        \caption{Euclidean distance between optimal values in the image space.}
        \label{subfig:const_disc_err_front}
    \end{subfigure}    
    \medskip    
    \begin{subfigure}[t]{.401\textwidth}
        \centering
        \includegraphics[width=\linewidth]{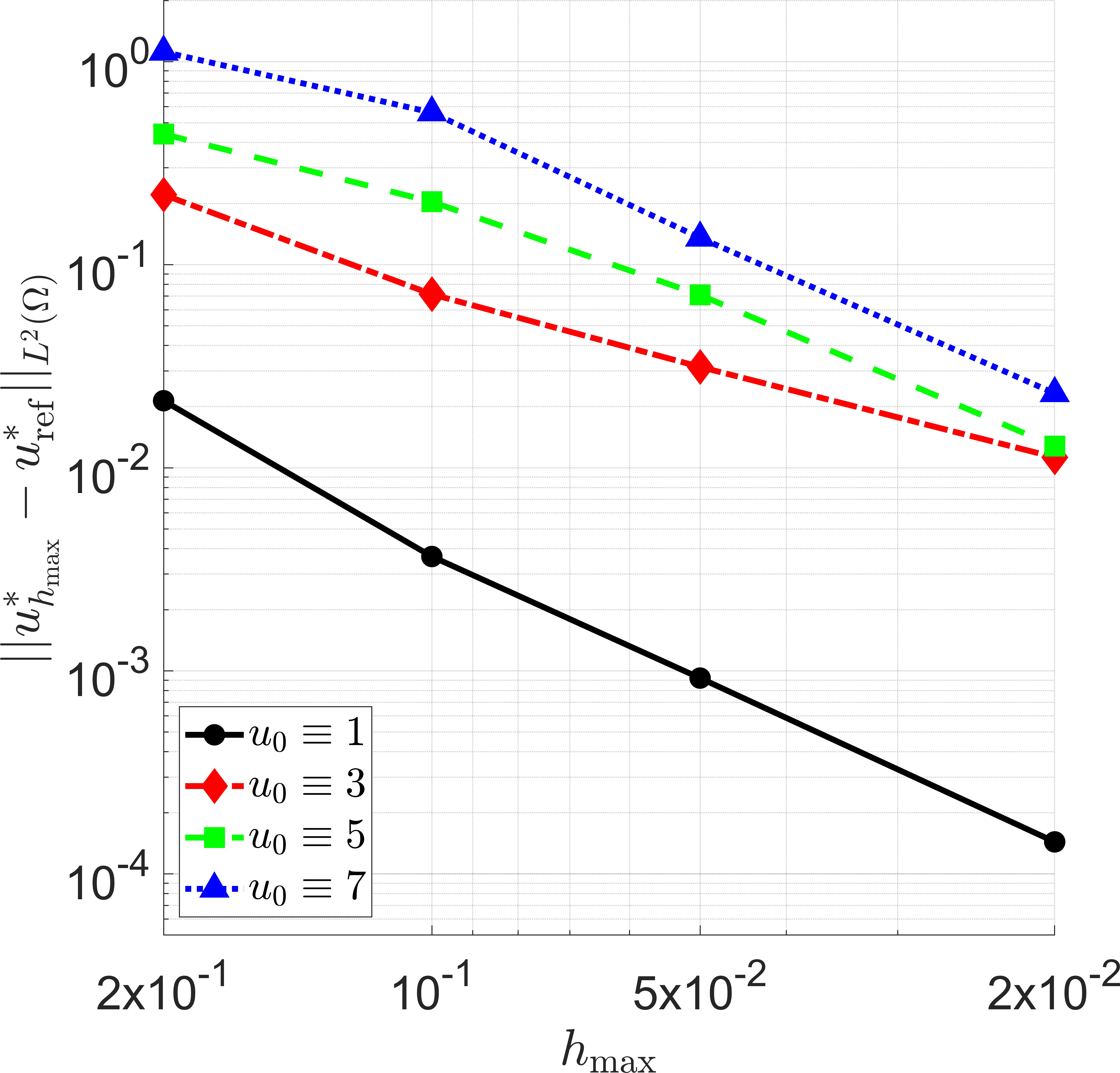}
        \caption{$L^2$-distance between optimal control and reference control.}
        \label{subfig:const_disc_err_control}
    \end{subfigure}
    \hspace{10mm}
    \begin{subfigure}[t]{.401\textwidth}
        \centering
        \includegraphics[width=\linewidth]{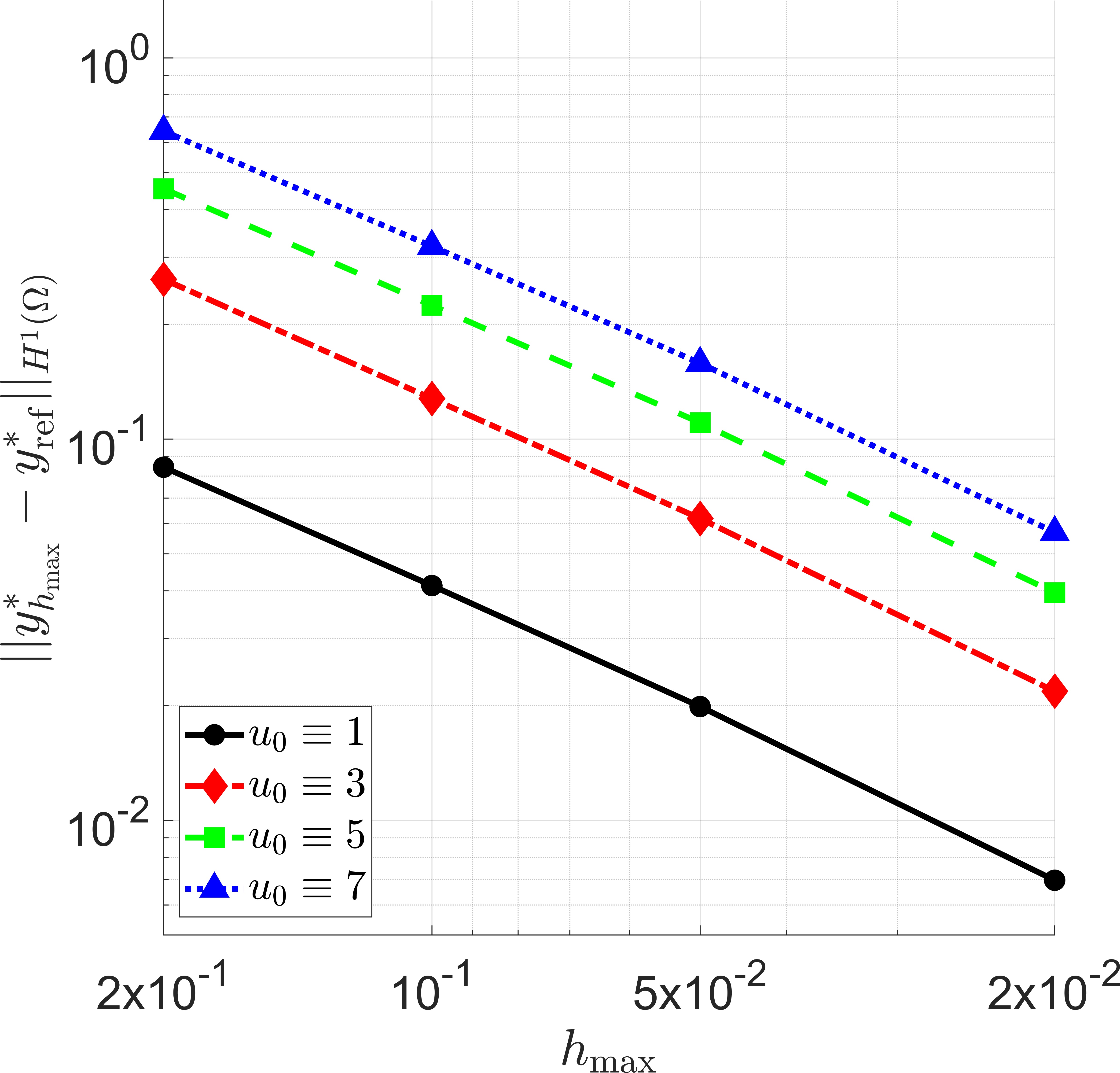}
        \caption{$H^1$-distance between optimal state and reference state.}
        \label{subfig:const_disc_err_state}
    \end{subfigure}
    \caption{Qualitative analysis of the solutions derived by Algorithm \ref{algo:nonsmooth_descent_method} for different discretizations for the constant obstacle. Subfigures (b) - (d) use the reference solution $u_{\mathrm{ref}}^*$ corresponding to mesh size $h_{\max} = 0.01$.}
    \label{fig:const_disc_err}
\end{center}
\end{figure}

A complete picture of the solutions obtained by Algorithm~\ref{algo:nonsmooth_descent_method} and the convergence behaviour is depicted in Figure~\ref{fig:const_disc_err} and Table~\ref{tab:const_num_iteration}. A qualitative analysis of the solutions is included in Figure \ref{fig:const_disc_err}. The iteration numbers required for each run are summarized in Table~\ref{tab:const_num_iteration}. For all initial values and mesh sizes the algorithm successfully terminates before reaching the maximum number of $10\,000$ iterations and computes an $(\overline{\varepsilon}, \overline{\delta})$-critical point. Subfigure~\ref{subfig:const_parefronts} shows the obtained solutions in the objective space for all initial values ranging from $u_0 \equiv 0$ to $u_0 \equiv 8$ and for all mesh sizes $h_{\max} \in \{ 0.2, 0.1, 0.05, 0.02 \}$ marked with different symbols and colors, respectively. The solutions with the same initial value (but for different mesh discretizations) cluster, while solutions for different initial values are evenly distributed and form a curved front. The clustering behaviour in the objective space will be examined further in Subfigure~\ref{subfig:const_disc_err_front}. The figure shows the distance of the objective function values of the obtained solutions to the objective function values of the reference solution. The plot contains one line for the different initial values $u_0 \equiv \hat{u} \in \{1, 3, 5, 7\}$ and shows how the distance evolves for finer meshes. Linear decay of the distances in double logarithmic scale can be observed, suggesting convergence of the front for $h_{\max} \to 0$. Similar behaviour can be observed in Subfigures~\ref{subfig:const_disc_err_control} and \ref{subfig:const_disc_err_state}. Subfigure \ref{subfig:const_disc_err_control} shows how the distance of the obtained control to the reference control in the $L^2$-norm evolves for finer meshes. The distance of the corresponding states to the reference state in the $H^1$-norm can be seen in Subfigure~\ref{subfig:const_disc_err_state}. In both subfigures, we can observe linear decay in the double logarithmic scale, indicating convergence of the controls and states computed by Algorithm~\ref{algo:nonsmooth_descent_method} for finer mesh sizes.

\begin{table}
    \centering
    \begin{tabular}{!{\vrule width 1.2pt}c|c!{\vrule width 1.2pt}r|r|r|r|r!{\vrule width 1.2pt}}
        
        \noalign{\hrule height 1.2pt}
        
        \multicolumn{2}{!{\vrule width 1.2pt}c!{\vrule width 1.2pt}}{}     & \multicolumn{5}{c!{\vrule width 1.2pt}}{$h_{\max}$}  \\  \cline{3-7}
        \multicolumn{2}{!{\vrule width 1.2pt}c!{\vrule width 1.2pt}}{}     & $0.2$      & $0.1$      & $0.05$     & $0.02$     & $0.01$ \\
        \noalign{\hrule height 1.2pt}
        \multirow{8}{*}{$u_0$}    & $1$ & $4$ & $4$        & $3$        & $3$        & $3$ \\ \cline{2-7}
        & $2$ & $4$        & $5$        & $4$ & $4$        & $4$ \\ \cline{2-7}
        & $3$ & $21$       & $22$       & $11$       & $37$       & $32$ \\ \cline{2-7}
        & $4$ & $106$      & $160$      & $236$      & $264$      & $91$ \\ \cline{2-7}
        & $5$ & $426$      & $2\,281$ & $1\,116$ & $542$      & $210$ \\ \cline{2-7}          
        & $6$ & $3\,885$ & $3\,894$ & $1\,619$ & $1\,025$ & $323$ \\ \cline{2-7}                
        & $7$ & $4\,756$ & $6\,190$ & $2\,918$ & $1\,370$ & $657$ \\ \cline{2-7}
        & $8$ & $2\,491$ & $3\,576$ & $3\,194$ & $2\,697$ & $822$ \\ \noalign{\hrule height 1.2pt}
    \end{tabular}
    \caption{Configuration 1: Number of iterations of Algorithm \ref{algo:nonsmooth_descent_method} for different initial values $u_0$ and mesh sizes $h_{\max}$.}
    \label{tab:const_num_iteration}
\end{table}

Table~\ref{tab:const_num_iteration} contains the number of iterations Algorithm~\ref{algo:nonsmooth_descent_method} performed for the different initial values and mesh sizes. For all mesh sizes the number of iterations increase with the magnitude of the initial control $u_0$. For $u_0 \equiv 1$ and $u_0 \equiv 2$ there is no contact between the state and the obstacle over the course of the optimization resulting in a small number of iterations. The number of iterations does not increase for finer meshes and we expect to converge to a finite value for $h_{\max} \to 0$ for all initial values.

\subsubsection{Configuration 2: Piecewise constant obstacle}
\label{subsub: Configuration 2 : piecewise constant obstacle}
        
\begin{figure}  
\begin{center}
\begin{subfigure}[t]{.28\textwidth}
        \centering
        \includegraphics[width=\linewidth]{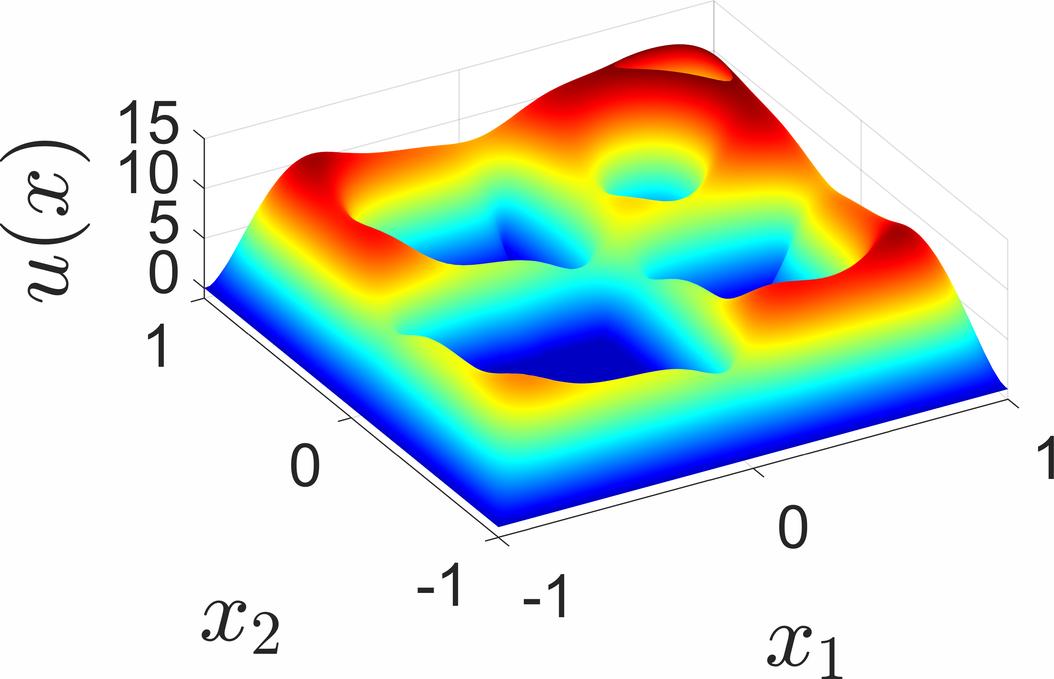}
        \caption{Control}
        \label{subfig:nonconst_control}
    \end{subfigure}
    \hfill
    \begin{subfigure}[t]{.28\textwidth}
        \centering
        \includegraphics[width=\linewidth]{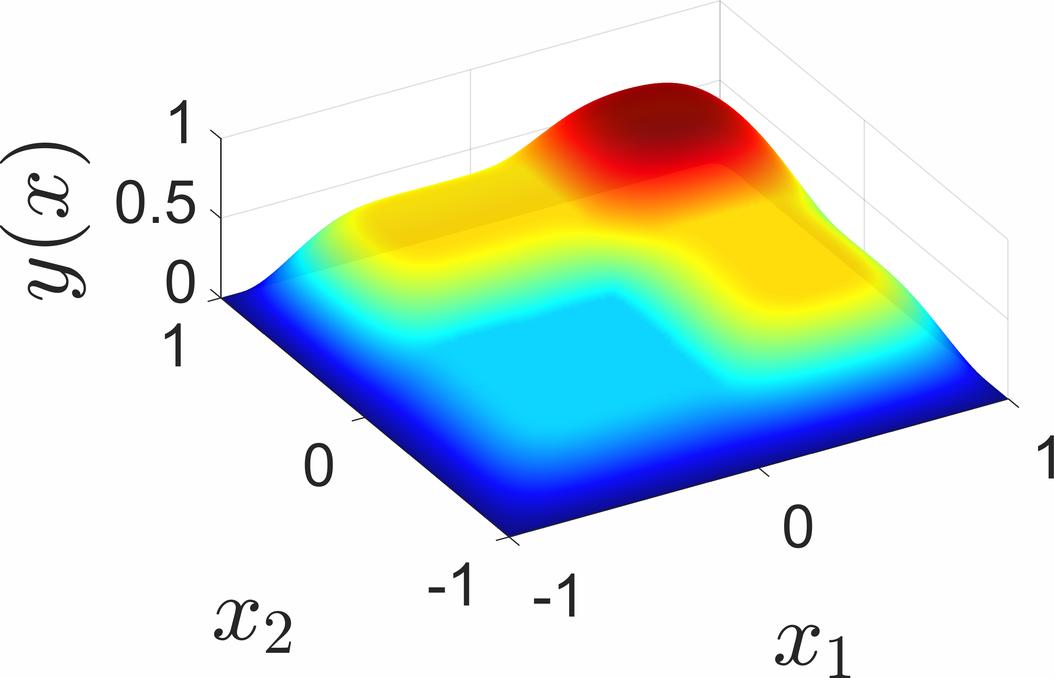}
        \caption{State}
        \label{subfig:nonconst_state}
    \end{subfigure}
    \hfill
    \begin{subfigure}[t]{.28\textwidth}
        \centering
        \includegraphics[width = \linewidth]{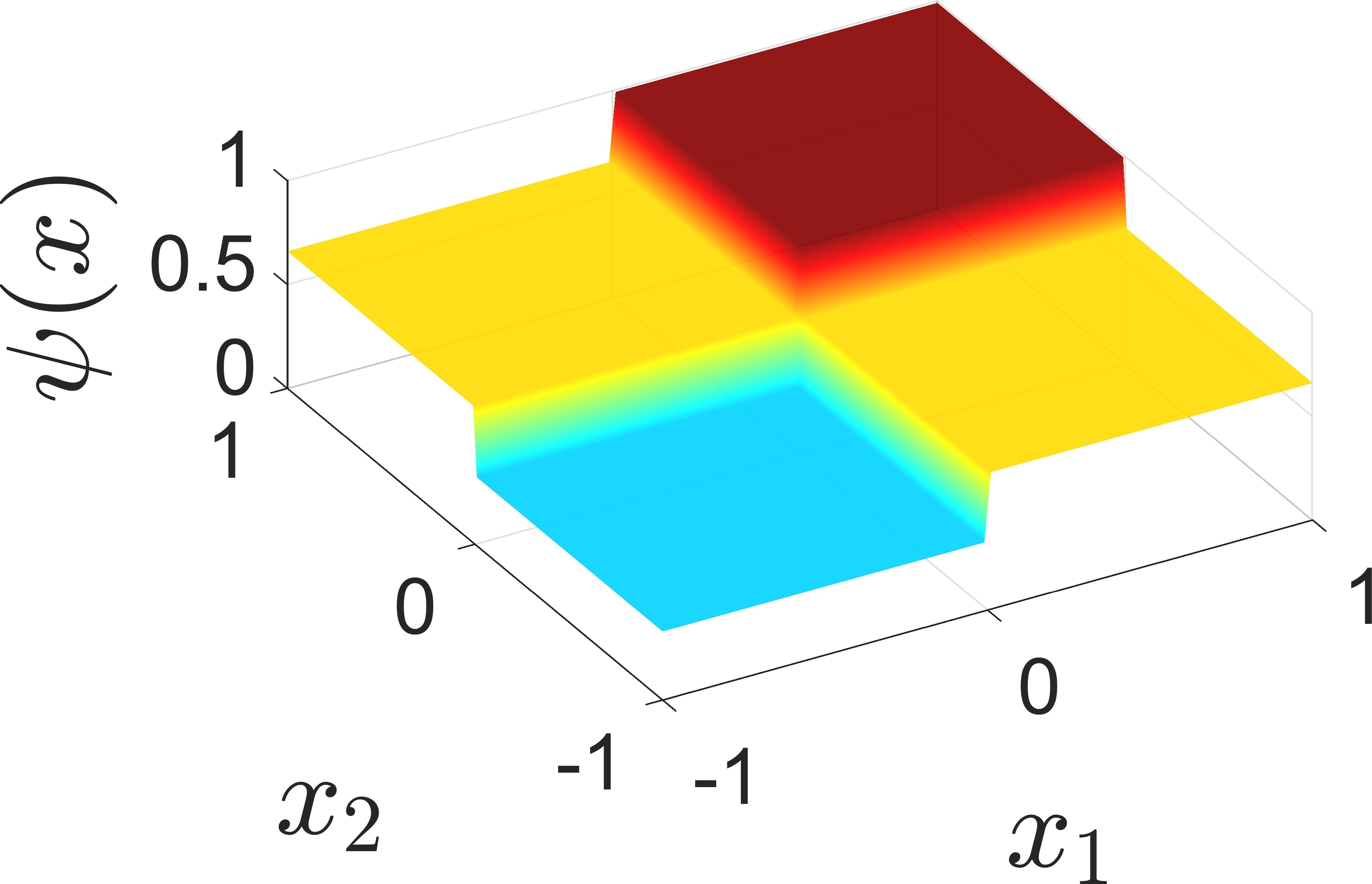}
        \caption{Obstacle}
        \label{subfig:nonconst_obstacle}
    \end{subfigure}
    \caption{A Pareto optimal control computed with Algorithm~\ref{algo:nonsmooth_descent_method} for mesh size $h_{\max} = 0.02$, initial control $u_0 \equiv 8$ and the piecewise constant obstalce $\psi$ defined in \eqref{eq:piecewise_constant_obstacle}.}
    \label{fig:nonconst_pareto_sol}    
\end{center}
\end{figure}

In the second example, we choose an obstacle $\psi$ given by a piecewise constant function defined by
\begin{align}
    \label{eq:piecewise_constant_obstacle}
    \psi \colon \Omega \to \R, \quad x \mapsto
    \begin{cases}
        \begin{aligned}
            &\nicefrac{1}{3} &&\text{if } x_1 \le 0 \text{ and } x_2 \le 0,\\
            & 1 &&\text{if } x_1 \ge 0 \text{ and } x_2 \ge 0,\\
            &\nicefrac{2}{3} &&\text{otherwise.}
        \end{aligned}
    \end{cases}
\end{align}
This obstacle can be interpreted analogously to that in Subsection~\ref{subsub: Configuration 1 : constant obstacle}. An approximate Pareto optimal control obtained by Algorithm~\ref{algo:nonsmooth_descent_method} for initial value $u_0 \equiv 8$ and $h_{\max} = 0.02$ together with the corresponding state can be seen in Subfigure~\ref{subfig:nonconst_control} and Subfigure~\ref{subfig:nonconst_state}. The obstacle $\psi$, defined in \eqref{eq:piecewise_constant_obstacle}, is shown in Subfigure~\ref{subfig:nonconst_obstacle}. Due to the nonconstant obstacle, we see a less structured behaviour in the control and state. Similarly to the first example, we observe vanishing control in areas with contact of the state with the obstacle.
    
Algorithmically, solving this problem configuration is expected to be more challenging compared to the first configuration with the constant obstacle, as the area of contact of the state changes more dynamically over the course of the algorithm's run, i.e., the problems nondifferentiability is more pronounced. 

\begin{figure}          
    \begin{center}    
    \label{fig:Discretization_errors_nonconstant_obstacle}
    \begin{subfigure}[t]{.401\textwidth}
        \centering
        \includegraphics[width = \linewidth]{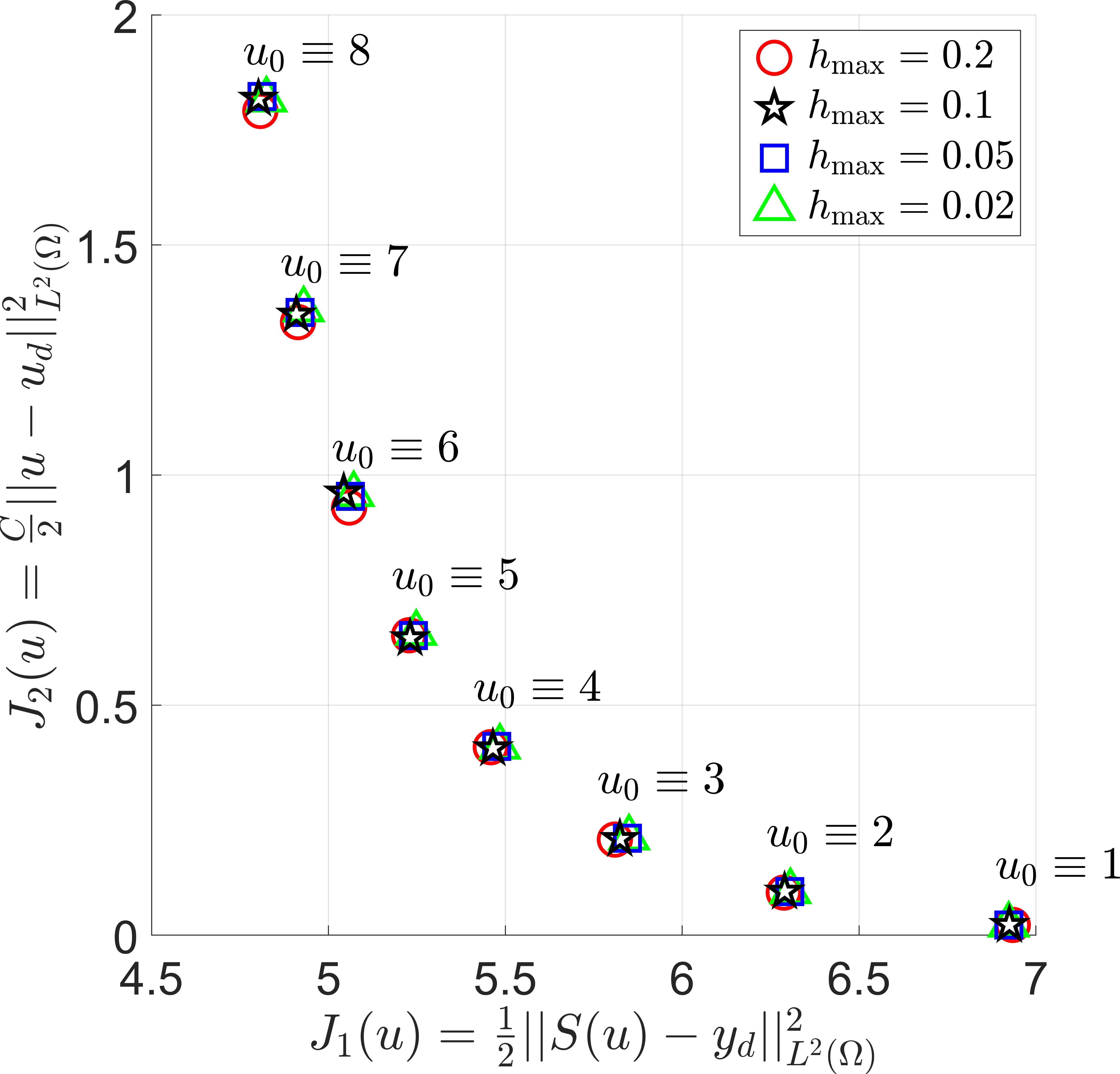}
        \caption{Pareto fronts for different mesh discretizations.}
        \label{subfig:nonconst_parefronts}
    \end{subfigure}
    \hspace{10mm}
    \begin{subfigure}[t]{.401\textwidth}
        \centering
        \includegraphics[width=\linewidth]{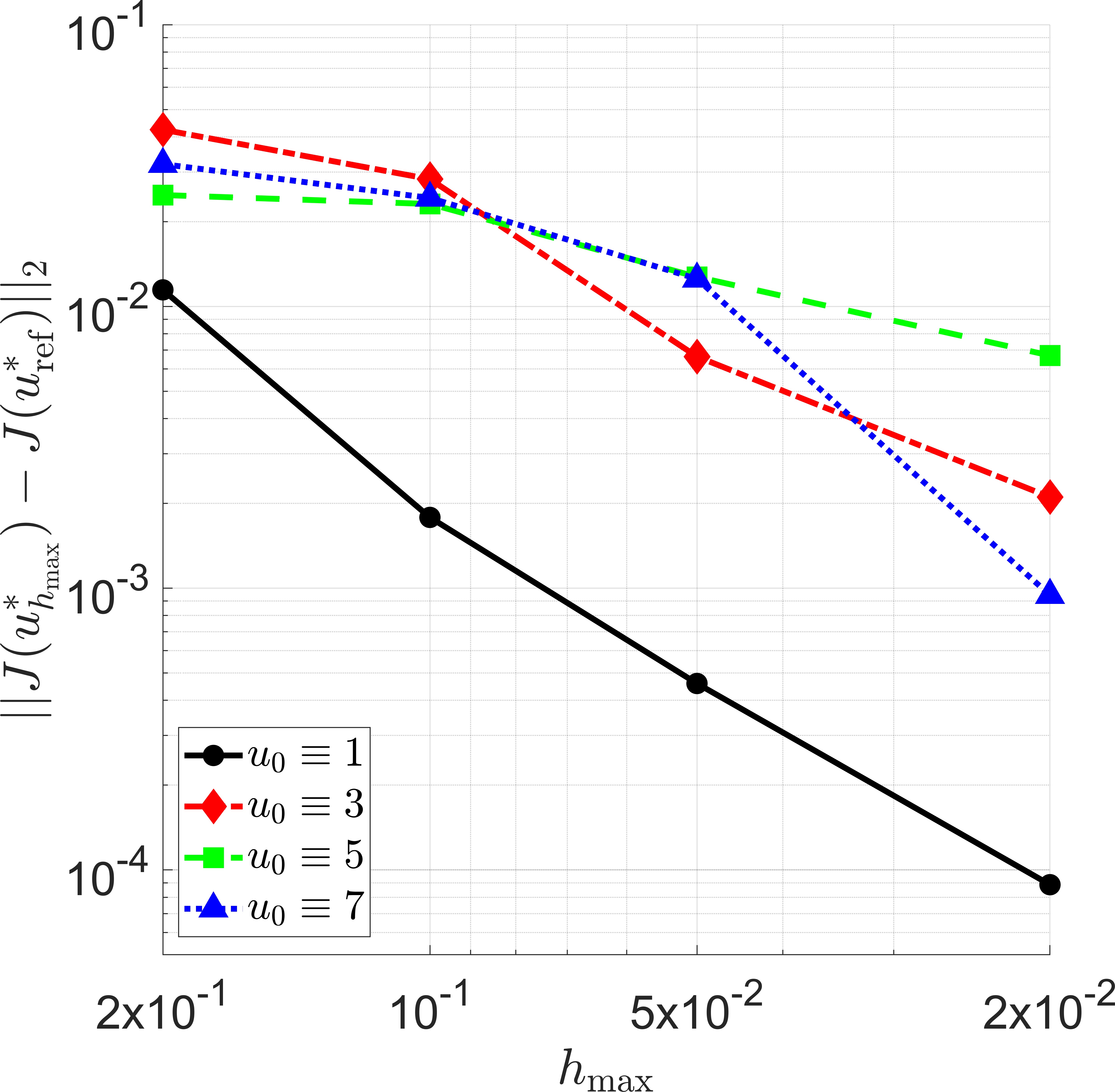}
        \caption{Euclidean distance of optimal values in the image space.}
        \label{subfig:nonconst_disc_err_front}
    \end{subfigure}    
    \medskip
    \begin{subfigure}[t]{.401\textwidth}
        \centering
        \includegraphics[width=\linewidth]{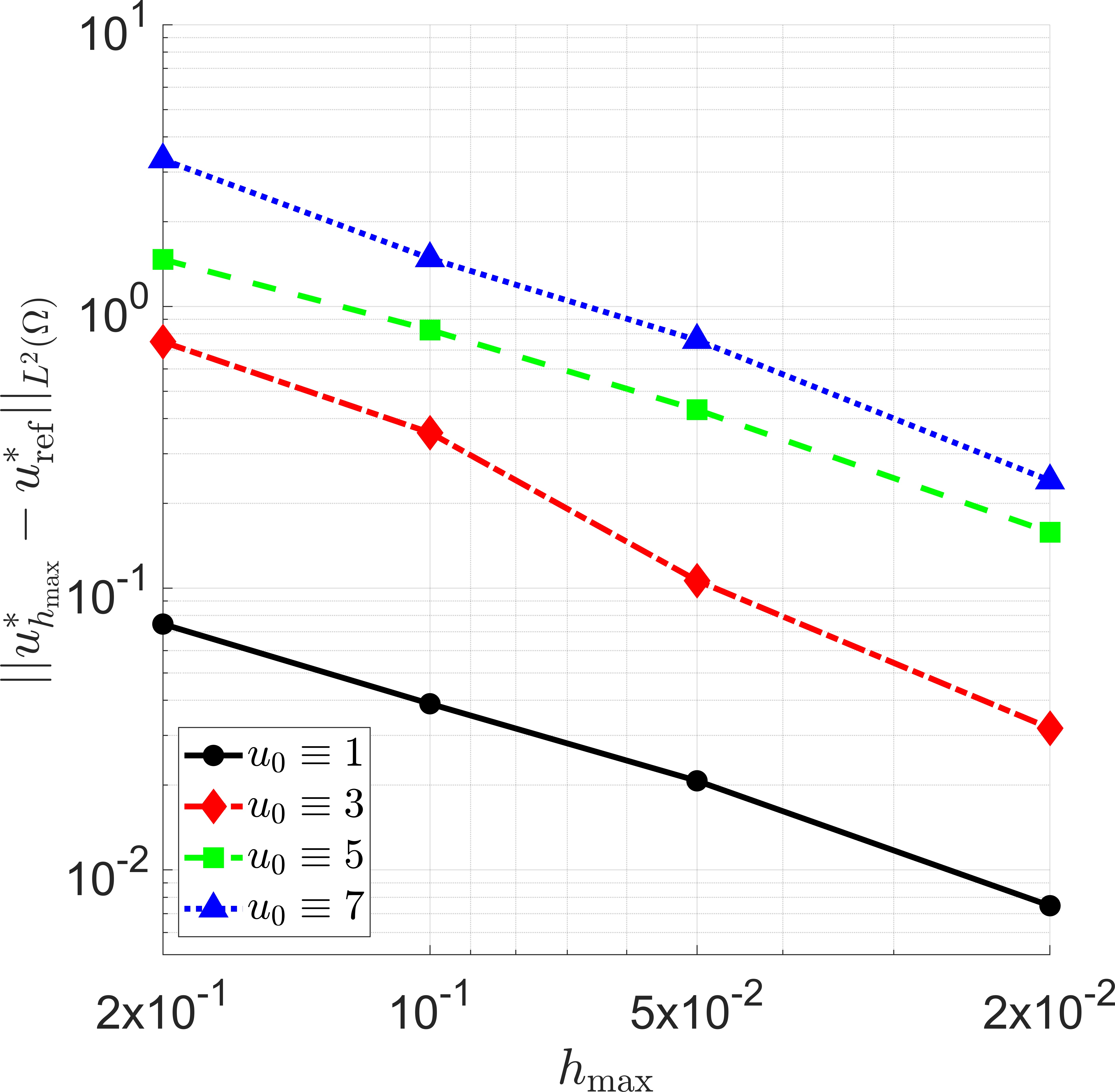}
        \caption{$L^2$-distance between optimal control and reference control.}
        \label{subfig:nonconst_disc_err_control}
    \end{subfigure}
    \hspace{10mm}
    \begin{subfigure}[t]{.401\textwidth}
        \centering
        \includegraphics[width=\linewidth]{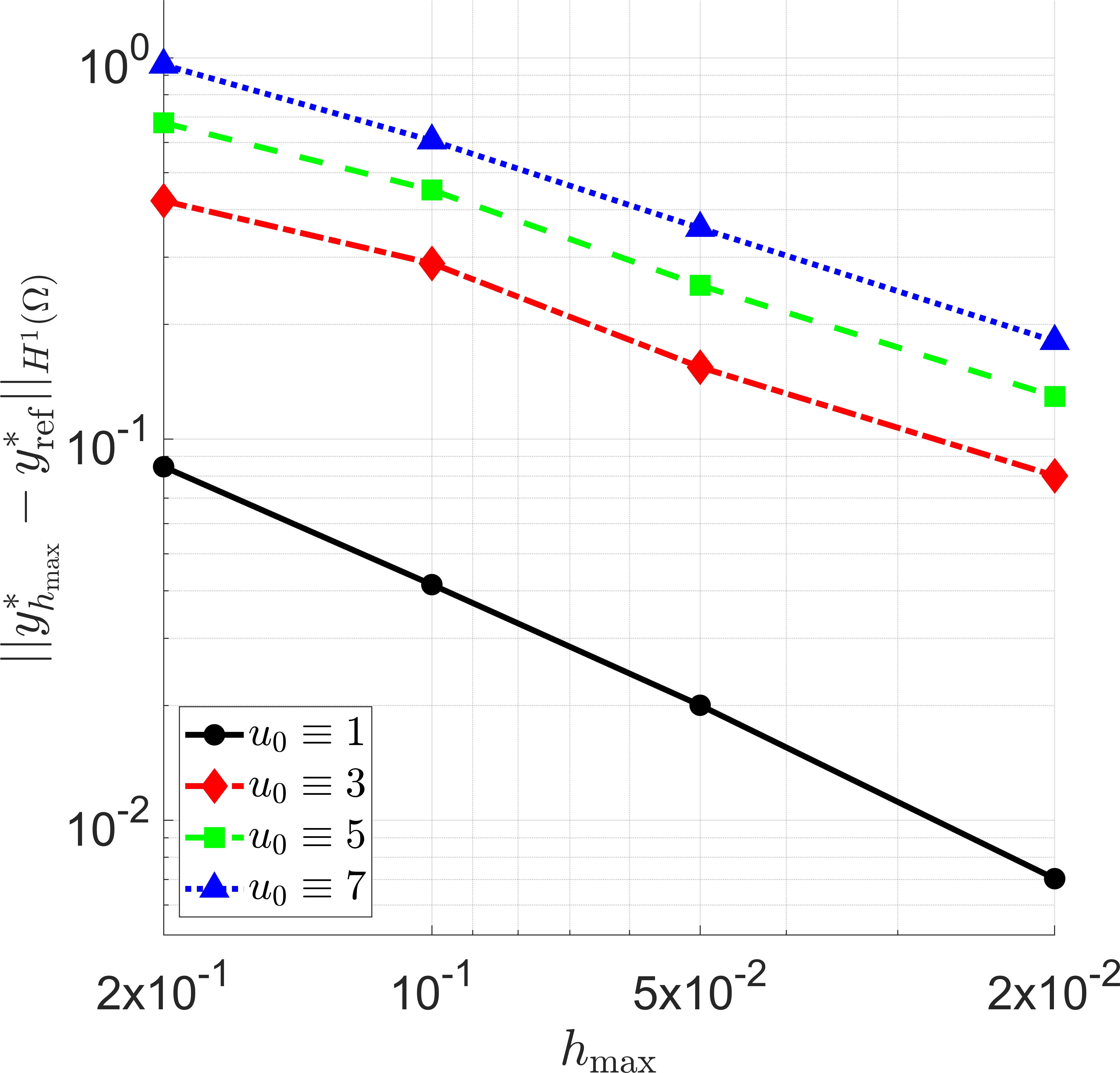}
        \caption{$H^1$-distance between optimal state and reference state.}
        \label{subfig:nonconst_disc_err_state}
    \end{subfigure}
    \caption{Qualitative analysis of the solutions derived by Algorithm~\ref{algo:nonsmooth_descent_method} for different discretizations for the nonconstant obstacle. Subfigures (b)-(d) use the reference solution $u_{\mathrm{ref}}^*$ corresponding to mesh size $h_{\max} = 0.01$.}
    \label{fig:nonconst_disc_err}
    \end{center}
\end{figure}

Figure~\ref{fig:nonconst_disc_err} contains a qualitative analysis of the solutions obtained by Algorithm~\ref{algo:nonsmooth_descent_method} for the piecewise constant obstacle. The objective function values obtained from Algorithm~\ref{algo:descent_direction} for different initial values $u_0 \equiv \hat u \in \{1, 2, \dots, 8\}$ and different mash sizes $h_{\max} \in \{0.2, 0.1, 0.05, 0.02 \}$ is visualized in Subfigure~\ref{subfig:nonconst_parefronts}. The objective functions form a front in the image space and solutions for different mesh discretizations but with same initial control cluster. This clustering is further examined in Subfigure \ref{subfig:nonconst_parefronts}, where the diminishing mesh size $h_{\max}$ is plotted over the distance between the computed objective function value and the objective function value of the reference solution. We observe linear decay of the distance in double logarithmic scale. Subfigures~\ref{subfig:nonconst_disc_err_control} and \ref{subfig:nonconst_disc_err_state} contain the distance of the obtained optimal control to reference control in the $L^2$-norm and the distance of the corresponding state to the reference stated in the $H^1$-norm, respectively. Again, we note linear decay for distances for smaller values of $h_{\max}$ in the double logarithmic scale. These plots indicate convergence of the solutions obtained by Algorithm~\ref{algo:nonsmooth_descent_method} for finer meshes.

Table~\ref{tab:nonconst_num_iteration} contains a comparison of the number of iterations performed to reach the stopping criterion in Algorithm~\ref{algo:nonsmooth_descent_method} for the different initial controls and the different meshes. We see the same trend as in the first example. However, for the piecewise constant obstacle, the iteration numbers are higher for almost all runs compared to the results for the constant obstacle. For all meshes we see an increasing number of iterations with an increasing magnitude of the initial control $u_0$. This is expected since for a higher magnitude of the initial control, we have more points with contact in the beginning. The number of iterations is bounded for the different mesh sizes and we expect convergence for $h_{\max} \to 0$ for all initial values of $u_0$. 

\begin{table}[]
    \centering
    \begin{tabular}{!{\vrule width 1.2pt}c|c!{\vrule width 1.2pt}r|r|r|r|r!{\vrule width 1.2pt}}
        \noalign{\hrule height 1.2pt}
        \multicolumn{2}{!{\vrule width 1.2pt}c!{\vrule width 1.2pt}}{}     & \multicolumn{5}{c!{\vrule width 1.2pt}}{$h_{\max}$}  \\  \cline{3-7}
        \multicolumn{2}{!{\vrule width 1.2pt}c!{\vrule width 1.2pt}}{}     & $0.2$      & $0.1$      & $0.05$     & $0.02$     & $0.01$ \\ \noalign{\hrule height 1.2pt}
        \multirow{8}{*}{$u_0$}    
            
        & $1$ & $10$     & $8$      & $37$     & $18$     & $17$ \\ \cline{2-7}
        & $2$ & $8$      & $8$      & $16$     & $22$     & $15$ \\ \cline{2-7}
        & $3$ & $41$     & $85$     & $48$     & $48$     & $28$ \\ \cline{2-7}
        & $4$ & $33$     & $863$    & $528$    & $286$    & $150$ \\ \cline{2-7}
        & $5$ & $324$    & $3\,135$ & $2\,381$ & $1\,254$ & $1\,013$ \\ \cline{2-7}
        & $6$ & $4\,070$ & $3\,701$ & $2\,696$ & $1\,046$ & $513$ \\ \cline{2-7}
        & $7$ & $9\,344$ & $7\,907$ & $5\,079$ & $1\,705$ & $827$ \\ \cline{2-7}
        & $8$ & $9\,719$ & $4\,539$ & $4\,757$ & $2\,387$ & $970$ \\ \noalign{\hrule height 1.2pt}
    \end{tabular}
    \caption{Configuration 2: Number of iterations of Algorithm~\ref{algo:nonsmooth_descent_method} for nonconstant obstacle.}
    \label{tab:nonconst_num_iteration}
\end{table}
    \clearpage
\subsubsection{Size of the approximated Goldstein $\varepsilon$-subdifferential}
\label{subsub:experiment_size_exp_subdiff}

\begin{figure}[h]
        \centering
        \includegraphics[width = .4\linewidth]{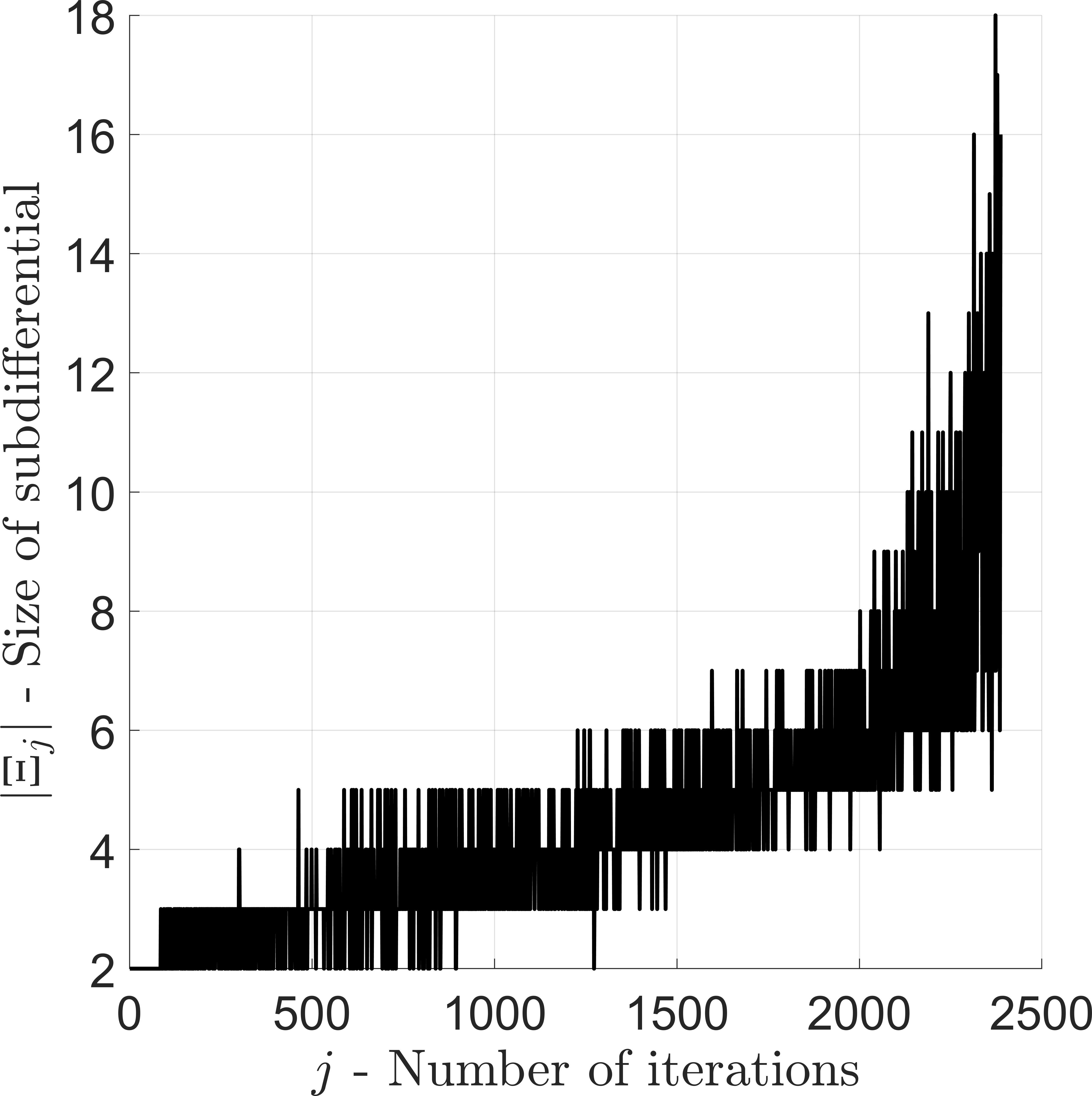}
        \caption{Size of the approximated subdifferential for each iteration. Results obtained by Algorithm \ref{algo:nonsmooth_descent_method} for the piecewise constant obstacle with mesh size $h_{\max} = 0.02$ and initial control $u_0 \equiv 8$.}
        \label{fig:Size_subdifferential}
    \end{figure}

In this subsection, we take a closer look at Step 2 in Algorithm~\ref{algo:nonsmooth_descent_method}. In this step a common descent direction yielding sufficient descent for all objective functions is computed using Algorithm~\ref{algo:descent_direction}. Algorithm~\ref{algo:descent_direction} computes a descent direction by iteratively updating an approximation $\Xi_l$ to the multiobjective $\varepsilon$-subdifferential, using subderivatives of the objective functions. Figure~\ref{fig:Size_subdifferential} shows the number of subderivatives in the final approximated $\varepsilon$-subdifferential in each iteration of a run of Algorithm~\ref{algo:nonsmooth_descent_method} with initial control $u_0 \equiv 8$ and maximum edge length $h_{\max} = 0.02$. We observe an increasing trend for the size of the subdifferential with the number of iterations. Up to iteration $900$ the algorithm regularly only requires two subderivatives. From iteration $1\,500$ onwards at least four subderivatives get used in every iteration. In the end, the subdifferential consists of up to $18$ subderivatives. This behaviour is not surprising: We expect the first objective function to be nonsmooth close to optima of the multiobjective control problem \eqref{prob:obstacle_red} (for the chosen initial control $u_0$), and hence, the algorithm converges to points, where the first objective function is not differentiable. To find a common descent direction in these areas, we need a sufficient number of subderivatives to describe the local behaviour of the objective function. The behaviour in Figure~\ref{fig:Size_subdifferential} can be observed across different mesh sizes and initial values. This indicates that the concept of Algorithm~\ref{algo:nonsmooth_descent_method} and the approximation of the multiobjective $\varepsilon$-subdifferential in Algorithm~\ref{algo:descent_direction} behave as expected.

	\section{Conclusion}
    \label{sec:conclusion_and_outlook}
    
    We extended a subderivative sampling approach to the setting of infinite-dimensional Hilbert spaces. To prove convergence of the presented algorithm to critical points, we extended the Goldstein $\varepsilon$-subdifferential to the multiobjective, infinite-dimensional setting. Its properties align with the ones known for the singleobjective and finite-dimensional cases. The main theorem on the multiobjective Goldstein $\varepsilon$-subdifferential describes a (set-valued) closedness property in the (strong $\times$ weak$^*$)-topology for a decaying sequence of scalars $\varepsilon_j \rightarrow 0$. This allows to conclude for a sequence $x_j \to x^*$ with $\xi_j \in F_{\varepsilon_j}(x_j)$ such that $\xi_j \rightharpoonup 0$ in the weak$^*$-topology that $x^*$ is in fact Pareto critical. In this environment we adapt the proof of convergence of the common descent method to the multiobjective, infinite-dimensional setting. The main results of this paper are Theorems~\ref{thm:convergence_proof_nonsmooth_descent_method} (which relies essentially on Theorem~\ref{thm:subdiff_eps_subdiff_multi}) and \ref{thm:convergence_proof_nonsmooth_descent_method_2}, showing that the common descent method defined in Algorithm \ref{algo:nonsmooth_descent_method} is capable of finding Pareto critical points and $(\varepsilon, \delta)$-critical points depending on the choice of algorithm parameters. To emphasize the usability of this algorithm, we apply it to a test problem configuration for multiobjective optimal control problems. The problem under consideration is formulated over an infinite-dimensional Hilbert space featuring a nondifferentiable objective function. Using a finite element approach we apply the nonsmooth descent method to the problem for different mesh sizes and initial controls. For the approximate Pareto front in the objective space, the control and the state, we observe mesh independent behaviour and convergence of the solutions. The number of iterations required for the algorithm does not grow for finer meshes. Furthermore, we have investigated the approximated multiobjective $\varepsilon$-subdifferentials used in the algorithm, showing that the concept for the approximation of the subdifferential works as expected.\\


\bibliographystyle{unsrt}
{\footnotesize
\bibliography{literature}
}

\end{document}